# THE FIRST PASSAGE EVENT FOR SUMS OF DEPENDENT LÉVY PROCESSES WITH APPLICATIONS TO INSURANCE RISK

By Irmingard Eder[1] and Claudia Klüppelberg

*Technische Universität München*

For the sum process $X = X^1 + X^2$ of a bivariate Lévy process $(X^1, X^2)$ with possibly dependent components, we derive a quintuple law describing the first upwards passage event of $X$ over a fixed barrier, caused by a jump, by the joint distribution of five quantities: the time relative to the time of the previous maximum, the time of the previous maximum, the overshoot, the undershoot and the undershoot of the previous maximum. The dependence between the jumps of $X^1$ and $X^2$ is modeled by a Lévy copula. We calculate these quantities for some examples, where we pay particular attention to the influence of the dependence structure. We apply our findings to the ruin event of an insurance risk process.

**1. Introduction.** In recent years, Lévy processes have been used to model multivariate insurance risk and operational risk processes, where the dependence between different business lines and risk types is crucial. The sum of the components of such a risk portfolio describes the total risk of an insurance company or a bank, and of special interest is the first passage event of the total risk process over a constant barrier, cf. [12, 18] for the univariate case and [4, 5, 13] for the multivariate case.

Motivated by these examples, we study in this paper the fluctuations of a one-dimensional (càdlàg) Lévy process $X = (X_t)_{t \geq 0}$, which is the sum of the components of a general multivariate Lévy process $\mathbf{X}$ in $\mathbb{R}^d$ with possibly dependent components. More precisely, we derive the joint distribution of those five quantities, which characterise the first upwards passage of $X$ over a constant barrier, when it is caused by a jump, and investigate also, which

Received December 2007; revised February 2009.
[1]Supported by the German Science Foundation (Deutsche Forschungsgemeinschaft) through the Graduiertenkolleg "Angewandte Algorithmische Mathematik."

component is likely to cause this passage. The paper generalizes results in [8, 13] for first passage events of Lévy processes.

We recall that the distribution of a $d$-dimensional Lévy process $\mathbf{X}$ is characterised by the Lévy–Khintchine representation

$$\mathbb{E}[e^{i(\mathbf{z},\mathbf{X}_t)}] = e^{-t\Psi(\mathbf{z})}, \qquad \mathbf{z} \in \mathbb{R}^d,$$

where

$$(1.1) \quad \Psi(\mathbf{z}) = i(\boldsymbol{\gamma},\mathbf{z}) + \frac{1}{2}\mathbf{z}^\top A\mathbf{z} + \int_{\mathbb{R}^d}(1 - e^{i(\mathbf{z},\mathbf{x})} + i(\mathbf{z},\mathbf{x})1_{\{|\mathbf{x}|\leq 1\}})\Pi(d\mathbf{x}),$$

$(\cdot,\cdot)$ denotes the inner product in $\mathbb{R}^d$, $\boldsymbol{\gamma} \in \mathbb{R}^d$ and $A$ is a symmetric nonnegative definite $d \times d$ matrix. The so-called *Lévy measure* $\Pi$ is a measure on $\mathbb{R}^d$ satisfying $\Pi((0,\ldots,0)) = 0$ and $\int_{\mathbb{R}^d} \min\{1, |\mathbf{x}|^2\}\Pi(d\mathbf{x}) < \infty$. A particular role is played by *subordinators*, which are Lévy processes, whose components have nondecreasing paths. Other important classes are *spectrally one-sided Lévy processes*, which have only positive or only negative jumps.

The important question, when considering multivariate Lévy processes, is the dependence modeling between the jumps of the components. Starting with the random walk model, a natural way of dependence modeling is by a copula, coupling the increments of the marginal random walks; see [14, 20]. This approach works also for compound Poisson processes (CPPes), but not for general Lévy processes with infinite Lévy measure. Consequently, we invoke a Lévy copula, which introduces dependence among the jump components of a multivariate Lévy process; see [6, 15].

We derive new results in fluctuation theory for sums of Lévy processes and study, in particular, the influence of jump dependence; for an introduction to fluctuation theory for Lévy processes we refer to [1], Chapter VI, and [19], Chapters 6 and 7. For an illustration, see Figure 7.1 of [19].

We shall need the following definitions, where all stochastic quantities in this paper are defined on a filtered probability space $(\Omega, \mathcal{F}, (\mathcal{F}_t)_{t\geq 0}, \mathbb{P})$.

We define the *running suprema* and *running infima of $X$* for $t > 0$

$$(1.2) \qquad \overline{X}_t := \sup_{u \leq t} X_u \quad \text{and} \quad \underline{X}_t := \inf_{u \leq t} X_u,$$

and the *first upwards passage time* over and the *first downwards passage time* under a fixed barrier $x \in \mathbb{R}$ by

$$(1.3) \qquad \tau_x^+ := \inf\{t > 0 : X_t > x\} \quad \text{and} \quad \tau_x^- := \inf\{t > 0 : X_t < x\}.$$

Further, we define the *time of the previous maximum of $X$* and the *time of the previous minimum of $X$* before time $t > 0$

$$(1.4) \quad \overline{G}_t := \sup\{s < t : \overline{X}_s = X_s\} \quad \text{and} \quad \underline{G}_t := \sup\{s < t : \underline{X}_s = X_s\}.$$

More precisely, we investigate the following quantities for a sum $X$ of possibly dependent Lévy processes, which characterize first upwards passage of $X$ over a fixed barrier caused by a jump:

(1) $\tau_x^+ - \overline{G}_{\tau_x^+-}$ time of first passage relative to time of previous maximum,
(2) $\overline{G}_{\tau_x^+-}$ time of previous maximum,
(3) $X_{\tau_x^+} - x$ overshoot,
(4) $x - X_{\tau_x^+-}$ undershoot, and
(5) $x - \overline{X}_{\tau_x^+-}$ undershoot of previous maximum.

Our paper is organized as follows. In Section 2, we consider the quintuple law in the random walk case as motivation. The general quintuple law for sums of possibly dependent Lévy processes is given in Section 3. In Section 4, two situations are considered, where all quantities of the quintuple law can be identified concretely. We calculate explicit quantities in Section 5, and give examples for different dependence structures. In Section 6, we apply our results to describe the ruin event in an insurance risk portfolio. The technical proofs are postponed to Section 7. For presentation purposes, we restrict ourselves to $d = 2$.

**2. The quintuple law for the sum of a bivariate random walk.** To see which jump of the Lévy process $(X^1, X^2)$ entails the first passage of the sum $X = X^1 + X^2$ and where the dependence affects this event, we decompose the jumps of $(X^1, X^2)$ in single, common, positive and negative jumps. We formulate first the quintuple law for the sum of a bivariate random walk $(Z_n^1, Z_n^2)_{n \in \mathbb{N}_0}$ starting in $(Z_0^1, Z_0^2) = 0$ given by

$$Z_n^1 = \sum_{i=1}^n \xi_i^1 \quad \text{and} \quad Z_n^2 = \sum_{i=1}^n \xi_i^2, \qquad n \in \mathbb{N},$$

where $(\xi_n^1, \xi_n^2)_{n \in \mathbb{N}}$ are i.i.d. with bivariate distribution function (d.f.) $F_{12}$ and margins $F_1$ and $F_2$, respectively. We are interested in first upwards passage across $x \geq 0$ of their sum

$$Z_0 = 0 \quad \text{and} \quad Z_n = \sum_{i=1}^n (\xi_i^1 + \xi_i^2), \qquad n \in \mathbb{N},$$

where $(\xi_n^1 + \xi_n^2)_{n \in \mathbb{N}}$ are i.i.d. with d.f. $F$. For $i = 1, 2$, we allow $F_i$ to have an atom at zero with the consequence that the random walks can have jumps of size 0 and so one marginal random walk can jump without the other. We separate the jumps of $Z$ according to their origin and their sign and decompose $Z$ for each $n \in \mathbb{N}$ into components as follows:

(2.1)
$$Z_n = P_n^1 + P_n^2 + P_n^3 + P_n^4 + P_n^5 + \sum_{i=1}^n \xi_i^1 1_{\{\xi_i^1 < 0, \xi_i^2 = 0\}}$$
$$+ \sum_{j=1}^n \xi_j^2 1_{\{\xi_j^1 = 0, \xi_j^2 < 0\}} + \sum_{k=1}^n (\xi_k^1 + \xi_k^2) 1_{\{\xi_k^1 < 0, \xi_k^2 < 0\}},$$

where $P^1, \ldots, P^5$ are those components, where upwards passage can happen:

$$P_n^1 = \sum_{i=1}^n \xi_i^1 1_{\{\xi_i^1 > 0, \xi_i^2 = 0\}}, \qquad P_n^2 = \sum_{i=1}^n \xi_i^2 1_{\{\xi_i^1 = 0, \xi_i^2 > 0\}},$$

$$P_n^3 = \sum_{i=1}^n (\xi_i^1 + \xi_i^2) 1_{\{\xi_i^1 > 0, \xi_i^2 > 0\}}, \qquad P_n^4 = \sum_{i=1}^n (\xi_i^1 + \xi_i^2) 1_{\{\xi_i^1 > 0, \xi_i^2 < 0\}},$$

$$P_n^5 = \sum_{i=1}^n (\xi_i^1 + \xi_i^2) 1_{\{\xi_i^1 < 0, \xi_i^2 > 0\}},$$

and the increments $\Delta P^k$ have d.f.s $F_{P^k}$ for $k = 1, \ldots, 5$. Further, we define the analogous quantities to (1.2)–(1.4): the *running maxima of $Z$* by

$$\overline{Z}_n := \max_{k \leq n} Z_k, \qquad n \in \mathbb{N}_0,$$

the *first strictly upwards passage time of $Z$* over a fixed barrier $x \in \mathbb{R}$

$$\mathcal{T}_x^+ := \min\{n \in \mathbb{N} : Z_n > x\},$$

and the *time of the previous maximum of $Z$* before time $n \in \mathbb{N}$

(2.2) $$\overline{G}^n := \max\{k \leq n : Z_k = \overline{Z}_n\}.$$

The quantities (1)–(5) from the Introduction are for the random walk $Z$:

(1) $\mathcal{T}_x^+ - 1 - \overline{G}^{\mathcal{T}_x^+ - 1}$ number of time points between previous maximum and first passage,

(2) $\overline{G}^{\mathcal{T}_x^+ - 1}$ time of previous maximum,

(3) $Z_{\mathcal{T}_x^+} - x$ overshoot,

(4) $x - Z_{\mathcal{T}_x^+ - 1}$ undershoot, and

(5) $x - \overline{Z}_{\mathcal{T}_x^+ - 1}$ undershoot of previous maximum.

Let $(T_n, H_n)_{n \in \mathbb{N}_0}$ be the weakly ascending and $(\widehat{T}_n^*, \widehat{H}_n^*)_{n \in \mathbb{N}_0}$ the strictly descending ladder process with *potential measures*

(2.3)
$$U(j, dx) = \sum_{n=0}^\infty \mathbb{P}(T_n = j, H_n \in dx),$$
$$\widehat{U}^*(i, dx) = \sum_{n=0}^\infty \mathbb{P}(\widehat{T}_n^* = i, \widehat{H}_n^* \in dx).$$

The proof of the following result is a consequence of the decomposition (2.1) in combination with the proof of Theorem 4 of [8].

THEOREM 2.1 (Quintuple law for the sum of a bivariate random walk). *Let $x > 0$ be a constant barrier. For $u > 0$, $y \in [0,x]$, $v \geq y$ and $i,j \in \mathbb{N}_0$ we have for $k = 1, \ldots, 5$,*

$$\mathbb{P}(\mathcal{T}_x^+ - 1 - \overline{G}^{\mathcal{T}_x^+ - 1} = i, \overline{G}^{\mathcal{T}_x^+ - 1} = j, Z_{\mathcal{T}_x^+} - x \in du, x - Z_{\mathcal{T}_x^+ - 1} \in dv,$$
$$x - \overline{Z}_{\mathcal{T}_x^+ - 1} \in dy, \Delta Z_{\mathcal{T}_x^+} = \Delta P_{\mathcal{T}_x^+}^k)$$
$$= F_{P^k}(du + v) \widehat{U}^*(i, dv - y) U(j, x - dy).$$

For the barrier $x = 0$, we have $\overline{Z}_{\mathcal{T}_0^+ - 1} = 0$ a.s. and get the following result.

COROLLARY 2.2. *Let $x = 0$ be a constant barrier. For $u > 0$, $v \geq 0$ and $i,j \in \mathbb{N}_0$, we have for $k = 1, \ldots, 5$,*

$$\mathbb{P}(\mathcal{T}_0^+ - 1 - \overline{G}^{\mathcal{T}_0^+ - 1} = i, \overline{G}^{\mathcal{T}_0^+ - 1} = j, Z_{\mathcal{T}_0^+} \in du,$$
$$-Z_{\mathcal{T}_0^+ - 1} \in dv, \Delta Z_{\mathcal{T}_0^+} = \Delta P_{\mathcal{T}_0^+}^k)$$
$$= F_{P^k}(du + v) \widehat{U}^*(i, dv) U(j, \{0\}).$$

REMARK 2.3. We can choose $\overline{G}_n^* = \min\{k \geq n : Z_k = \overline{Z}_n\}$ instead of $\overline{G}^n$ in (2.2). Then we work with the strictly ascending (instead of weakly ascending) and the weakly descending (instead of strictly descending) ladder processes.

We have not yet specified the dependence structure between the random walks $Z^1$ and $Z^2$. Since we want to study the influence of different dependence structures, we model dependence by some copula on the increments $\xi^1$ and $\xi^2$; see [14, 20] for details. By *Sklar's theorem* for copulas (cf. [20], Theorem 2.3.3), we write the joint d.f. of $(\xi^1, \xi^2)$ as

(2.4) $$F_{12}(x,y) = C(F_1(x), F_2(y)), \qquad x,y \in [-\infty, \infty],$$

where the copula $C$ is unique, if $F_1$ and $F_2$ are continuous; otherwise, $C$ is uniquely determined on $\mathrm{Ran}\, F_1 \times \mathrm{Ran}\, F_2$. Then we find expressions for the d.f. $F$ of the sum $\xi^1 + \xi^2$ and also of $F_{P^k}$, which makes the quintuple law of Theorem 2.1 and Corollary 2.2 precise in reference of the chosen copula. In the following result, we only consider the situation, where both random walks always jump together. If $F_1$, $F_2$ have atoms in 0, then we decompose the random walks as in (2.1) and observe that the absolutely continuous parts of $F_1$ and $F_2$ may have total mass smaller than 1.

THEOREM 2.4. *Suppose that $F_i$ for $i = 1,2$ are absolutely continuous and the dependence between $Z^1$ and $Z^2$ is modeled by a twice continuously* *differentiable copula $C$. Then $P^1 = P^2 = 0$ a.s. and $F_{P^k}$ for $k = 3, 4, 5$ of Theorem 2.1 are given for $z > 0$ by*

$$F_{P^3}(z) = \int_0^z \left[ \frac{\partial C(u,v)}{\partial u} \bigg|_{u=F_1(x_1)} \right]_{F_2(0)}^{F_2(z-x_1)} F_1(dx_1),$$

(2.5) $$F_{P^4}(z) = \int_{-\infty}^0 \left[ \frac{\partial C(u,v)}{\partial v} \bigg|_{v=F_2(x_2)} \right]_{F_1(0)}^{F_1(z-x_2)} F_2(dx_2),$$

$$F_{P^5}(z) = \int_{-\infty}^0 \left[ \frac{\partial C(u,v)}{\partial u} \bigg|_{u=F_1(x_1)} \right]_{F_2(0)}^{F_2(z-x_1)} F_1(dx_1).$$

PROOF. Since $F_1, F_2$ are absolutely continuous, all increments of $Z^1$ and $Z^2$ are nonzero and $P^1 = P^2 = 0$ a.s. From (2.4), we obtain for $x_1, x_2 \in \mathbb{R}$

$$F_{12}(dx_1, dx_2) = \frac{\partial^2 C(u,v)}{\partial u \, \partial v} \bigg|_{u=F_1(x_1), v=F_2(x_2)} F_1(dx_1) F_2(dx_2).$$

Furthermore, $F_{P^3}(z) = \int_0^z \int_0^{z-x_1} F_{12}(dx_1, dx_2)$ for $z > 0$. The expressions for $F_{P^4}$ and $F_{P^5}$ follow analogously. □

The potential measures $U$ and $\widehat{U}^*$ in Theorem 2.1 can be identified only in special cases.

Recall the $n$-fold convolution $F^{n*}(dx)$ of a probability measure $F(dx)$, where $F^{0*}(dx) = \delta_0(dx)$ is the Dirac-measure in 0 and $F^{1*} = F$.

THEOREM 2.5. *Suppose that for $i = 1, 2$ the random walks $Z^i$ have only positive jumps. Further, let the $F_i$ be absolutely continuous and the dependence between $Z^1$ and $Z^2$ is modeled by a twice continuously differentiable copula $C$. Then $\widehat{U}^*(\{0\}, \{0\}) = 1$ and, for $j \in \mathbb{N}_0$ and $x \geq 0$,*

$$U(j, dx) = F_{P^3}^{j*}(dx),$$

*where $F_{P^3}$ is given in (2.5).*

PROOF. $Z$ reaches a new maximum with every jump. So in (2.3), we have $\mathbb{P}(T_n = j, H_n \in dx) = 1_{\{n=j\}} \mathbb{P}(H_j \in dx)$ and $H_j$ is the sum of $j$ independent jumps with d.f. $F_{P^3}$. □

**3. The quintuple law for the sum of two Lévy processes.** For an arbitrary bivariate Lévy process $(X^1, X^2)$, we consider $X = X^1 + X^2$, which is again a Lévy process; see [21], Proposition 11.10. The proofs of our results rely on the Lévy–Itô decomposition of $(X^1, X^2)$ into independent parts, corresponding to (1.1),

$$\mathbf{X}_t = \mathbf{W}_t + \mathbf{S}_t, \qquad t \geq 0,$$

where $\mathbf{W}$ is a Gaussian process in $\mathbb{R}^2$ with characteristic triple $(\boldsymbol{\gamma}, A, 0)$. The Lévy process $\mathbf{S}$ is the jump part of $\mathbf{X}$ with Lévy measure $\Pi$ and has representation

$$
\begin{aligned}
\mathbf{S}_t = & \int_{(0,t]} \int_{|\mathbf{x}|>1} \mathbf{x} J(d\mathbf{x}, ds) \\
& + \lim_{\varepsilon \downarrow 0} \int_{(0,t]} \int_{\varepsilon<|\mathbf{x}|\leq 1} (\mathbf{x} J(d\mathbf{x}, ds) - \mathbf{x} \Pi(d\mathbf{x})\, ds), \qquad t \geq 0,
\end{aligned}
\tag{3.1}
$$

see [21], Theorem 19.2. The convergence in the second term on the right-hand side is a.s. and uniformly on compacts for $t \in [0, \infty)$. The measure $J$ is a Poisson random measure with intensity measure $\Pi(d\mathbf{x})\, ds$ on $\mathbb{R}^2 \times (0, \infty)$. We investigate first upwards passage by a jump of the sum process $X$, equivalently by a jump of $S = S^1 + S^2$.

Analogously to the random walk in Section 2, we decompose the paths of $(S^1, S^2)$ according to their jump behavior in single, common, positive and negative jumps. This causes no problem, if $(S^1, S^2)$ has a.s. sample paths of bounded variation; see [19], Exercise 2.8. But $(S^1, S^2)$ may have a.s. sample paths of unbounded variation, so that, according to [2] relation (31.32), a pathwise decomposition is not possible. In this case, we truncate for arbitrary $0 < \varepsilon < 1$ all jumps smaller than $\varepsilon$ and consider first the process for $t \geq 0$

$$
\begin{pmatrix} S_t^{1,\varepsilon} \\ S_t^{2,\varepsilon} \end{pmatrix} = \int_{(0,t]} \int_{|\mathbf{x}|>1} \mathbf{x} J(d\mathbf{x}, ds) \\
+ \int_{(0,t]} \int_{\varepsilon<|\mathbf{x}|\leq 1} (\mathbf{x} J(d\mathbf{x}, ds) - \mathbf{x} \Pi(d\mathbf{x})\, ds),
\tag{3.2}
$$

which is a CPP with drift $(D_{S^{1,\varepsilon}}, D_{S^{2,\varepsilon}}) = -\int_{\varepsilon<|\mathbf{x}|\leq 1} \mathbf{x}\Pi(d\mathbf{x})$ and Lévy measure $\Pi(\cdot) 1_{\{|\mathbf{x}|>\varepsilon\}}$. For $\varepsilon \downarrow 0$, the process $(S_t^{1,\varepsilon}, S_t^{2,\varepsilon})_{t\geq 0}$ converges a.s. to $(S_t^1, S_t^2)_{t\geq 0}$ and the convergence is locally uniformly in $t \in [0, \infty)$; see Lemma 20.7 of [21]. Now, we can decompose $(S^{1,\varepsilon}, S^{2,\varepsilon})$ in independent components. We denote by $S^{1,\varepsilon,+}$ the process of single positive jumps of $S^{1,\varepsilon}$; i.e., for $t > 0$

$$
S_t^{1,\varepsilon,+} = \int_{(0,t]} \int_{x_1>1} x_1 J((dx_1, \{0\}), ds) \\
+ \int_{(0,t]} \int_{\varepsilon<x_1\leq 1} (x_1 J((dx_1, \{0\}), ds) - x_1 \Pi(dx_1, \{0\})\, ds)
$$

and by $S^{1,\varepsilon,-}$ the single negative jumps of $S^{1,\varepsilon}$, i.e., for $t > 0$

$$
S_t^{1,\varepsilon,-} = \int_{(0,t]} \int_{x_1<-1} x_1 J((dx_1, \{0\}), ds) \\
+ \int_{(0,t]} \int_{-1\leq x_1<-\varepsilon} (x_1 J((dx_1, \{0\}), ds) - x_1 \Pi(dx_1, \{0\})\, ds).
$$

$S^{2,\varepsilon,+}$ and $S^{2,\varepsilon,-}$ are defined analogously for $S^{2,\varepsilon}$.

The processes $S^{1,\varepsilon,ij}$ and $S^{2,\varepsilon,ij}$ for $i,j \in \{+,-\}$ are the dependent jump parts of $(S^{1,\varepsilon}, S^{2,\varepsilon})$, where e.g., $S^{1,\varepsilon,++}$ denotes the positive jumps of $S^{1,\varepsilon}$, which happen together with positive jumps of $S^{2,\varepsilon}$; i.e., for $t > 0$,

$$S_t^{1,\varepsilon,++} = \int_{(0,t]} \int_{x_1 > 1} x_1 J((dx_1,(0,\infty)), ds)$$
$$+ \int_{(0,t]} \int_{\varepsilon < x_1 \leq 1} (x_1 J((dx_1,(0,\infty)), ds) - x_1 \Pi(dx_1,(0,\infty))\, ds).$$

Analogously, $S^{2,\varepsilon,++}$ denotes the positive jumps of $S^{2,\varepsilon}$, which happen together with positive jumps of $S^{1,\varepsilon}$. The notations $S^{1,\varepsilon,+-}$, $S^{2,\varepsilon,+-}$ and $S^{1,\varepsilon,--}$, $S^{2,\varepsilon,--}$ should be clear now.

This implies the following Lévy–Itô decomposition for the sum process

$$X = X^1 + X^2 = W^1 + W^2 + S^1 + S^2$$
$$(3.3) \quad = W^1 + W^2$$
$$+ \lim_{\varepsilon \downarrow 0} (P^{1,\varepsilon} + P^{2,\varepsilon} + P^{3,\varepsilon} + P^{4,\varepsilon} + P^{5,\varepsilon} + S^{1,\varepsilon,-} + S^{2,\varepsilon,-} + S^{\varepsilon,--}),$$

where $W^1 + W^2$ denotes the Gaussian part of $X$, which is independent of the jump component, and in (3.3) we have set $S^{\varepsilon,--} := S^{1,\varepsilon,--} + S^{2,\varepsilon,--}$. Then we summarize, using analogous notation to (2.1):

$$P^{1,\varepsilon} := S^{1,\varepsilon,+}, \qquad P^{2,\varepsilon} := S^{2,\varepsilon,+}, \qquad P^{3,\varepsilon} := S^{1,\varepsilon,++} + S^{2,\varepsilon,++},$$
$$P^{4,\varepsilon} := S^{1,\varepsilon,+-} + S^{2,\varepsilon,+-}, \qquad P^{5,\varepsilon} := S^{1,\varepsilon,-+} + S^{2,\varepsilon,-+},$$

which are all independent Lévy processes, because they a.s. never jump together; see [21], Exercise 12.10 on page 67. Since all processes in (3.3) are independent, we can let $\varepsilon \downarrow 0$ componentwise. According to [21], Lemma 20.7, we have

$$(3.4) \qquad \lim_{\varepsilon \downarrow 0} P^{k,\varepsilon} =: P^k \qquad \text{a.s.,}$$

where the convergence is uniformly on compacts for $t \in [0,\infty)$ for $k = 1,\ldots,5$.

Our main result is derived as a consequence of the Wiener–Hopf factorization, which is based on ladder processes. Using the same notation as in [8, 19], we denote by $(L_t)_{t \geq 0}$ and $(\widehat{L}_t)_{t \geq 0}$ the *local time at the maximum* and *at the minimum* and by $(L_t^{-1}, H_t)_{t \geq 0}$ and $(\widehat{L}_t^{-1}, \widehat{H}_t)_{t \geq 0}$ the (possibly killed) bivariate subordinators representing the *ascending* and *descending ladder processes*. Recall from [1], Proposition VI.4, and [19], page 158, that with the exception of a CPP all local extrema of $X$ are distinct. Therefore, we exclude CPPes in the following and treat them separately in Remark 3.4. The following situation is for every Lévy process $X$, which is not a CPP.

The joint Laplace exponents $\kappa$ and $\widehat{\kappa}$ of the above subordinators are for $\alpha, \beta \geq 0$ defined by the identities

$$e^{-\kappa(\alpha,\beta)} = \mathbb{E}[e^{-\alpha L_1^{-1} - \beta H_1} 1_{\{1 < L_\infty\}}],$$

$$e^{-\widehat{\kappa}(\alpha,\beta)} = \mathbb{E}[e^{-\alpha \widehat{L}_1^{-1} - \beta \widehat{H}_1} 1_{\{1 < \widehat{L}_\infty\}}],$$

where $L_\infty := \lim_{t \to \infty} L_t$ and $\widehat{L}_\infty := \lim_{t \to \infty} \widehat{L}_t$. By equations (6.15) and (6.16) of [19], we can also write for $\beta \in [0, \infty) + i\mathbb{R}$,

$$(3.5) \qquad \kappa(0, \beta) = q + \xi(\beta) = q + D_H \beta + \int_{(0,\infty)} (1 - e^{-\beta x}) \Pi_H(dx),$$

where $q \geq 0$ is the *killing rate* of $H$ so that $q > 0$ if and only if $\lim_{t \to \infty} X_t = -\infty$ a.s., $D_H = -\gamma_H - \int_{|x| \leq 1} x \Pi_H(dx) \geq 0$ is the drift of $H$, and $\Pi_H$ its Lévy measure. Note that the function $\xi(\cdot)$ is the Laplace exponent of an unkilled subordinator. Similar notation is used for $\widehat{\kappa}(0, \beta)$ by replacing $q$, $\xi$, $D_H$ and $\Pi_H$ by $\widehat{q}$, $\widehat{\xi}$, $D_{\widehat{H}}$ and $\Pi_{\widehat{H}}$. We also recall that whenever $q > 0$, we have $\widehat{q} = 0$.

Associated with the ascending and descending ladder processes are the *bivariate potential measures* on $[0, \infty)^2$

$$(3.6) \qquad \mathcal{U}(ds, dx) = \int_0^\infty \mathbb{P}(L_t^{-1} \in ds, H_t \in dx)\, dt,$$

$$(3.7) \qquad \widehat{\mathcal{U}}(ds, dx) = \int_0^\infty \mathbb{P}(\widehat{L}_t^{-1} \in ds, \widehat{H}_t \in dx)\, dt.$$

We recall that local time at the maximum is defined only up to a multiplicative constant; see e.g., [19], page 190, or [1], page 110. As a consequence, the exponent $\kappa$ can only be defined up to a multiplicative constant, which is then also inherited by $\mathcal{U}$.

Now, we are ready to state our first main result. Its proof can be found in Section 7. We recall the definition of the $P^k$ in (3.4).

THEOREM 3.1 (Quintuple law for the sum of Lévy processes). *Suppose that $X$ is not a CPP and $\Pi_{S^1}((0, \infty))$, $\Pi_{S^2}((0, \infty)) > 0$. Consider first upwards passage of $X$ over a constant barrier $x > 0$. Then there exists a normalization of local time at the maximum, given by the identity*

$$(3.8) \qquad q = \kappa(q, 0)\widehat{\kappa}(q, 0), \qquad q \geq 0,$$

*such that for $u > 0, y \in [0, x], v \geq y, s \geq 0, t \geq 0$, and for $k = 1, \ldots, 5$,*

$$(3.9) \quad \begin{aligned} \mathbb{P}(\tau_x^+ - \overline{G}_{\tau_x^+ -} &\in dt, \overline{G}_{\tau_x^+ -} \in ds, X_{\tau_x^+} - x \in du, x - X_{\tau_x^+ -} \in dv, \\ & x - \overline{X}_{\tau_x^+ -} \in dy, \Delta X_{\tau_x^+} = \Delta P_{\tau_x^+}^k) \\ &= \Pi_{P^k}(du + v)\widehat{\mathcal{U}}(dt, dv - y)\mathcal{U}(ds, x - dy). \end{aligned}$$

For the barrier $x = 0$, the situation simplifies by considering

(R): 0 is regular for $(0, \infty)$; i.e., $\tau_0^+ = 0$ a.s., or
(I): 0 is irregular for $(0, \infty)$; i.e., $\tau_0^+ > 0$ a.s.

Since we still exclude that $X$ is CPP, (R) holds if and only if (see [19], Theorem 6.5):

- $S^1 + S^2$ is of unbounded variation, or
- $S^1 + S^2$ is of bounded variation and
  – $X$ has a Gaussian component, or
  – $X$ has non-Gaussian component, but
    * $X$ has drift $D_X = -\gamma_X - \int_{|x| \le 1} x \Pi_X(dx) > 0$, or
    * $X$ has drift $D_X = 0$ and $\int_0^1 \frac{x}{\int_0^x \Pi_X((-\infty, -y)) \, dy} \Pi_X(dx) = \infty$.

(I) holds if and only if $S^1 + S^2$ is of bounded variation, $X$ has non-Gaussian component and either

- $D_X < 0$, or
- $D_X = 0$ and $\int_0^1 \frac{x}{\int_0^x \Pi_X((-\infty, -y)) \, dy} \Pi_X(dx) < \infty$.

COROLLARY 3.2. *Suppose that $X$ is not a CPP and $\Pi_{S^1}((0, \infty)), \Pi_{S^2}((0, \infty)) > 0$. Consider first passage of the barrier $x = 0$. Then $-\overline{X}_{\tau_0^+-} = 0$ a.s. and $\overline{G}_{\tau_0^+-} = 0$ a.s.*

(1) *If* (R) *holds, then*
$$-\overline{X}_{\tau_0^+-} = \overline{G}_{\tau_0^+-} = \tau_0^+ = X_{\tau_0^+} = -X_{\tau_0^+-} = 0 \qquad a.s.$$

(2) *If* (I) *holds, then for $u > 0, t \ge 0, v \ge 0$, and for $k = 1, \ldots, 5$,*
$$\mathbb{P}(\tau_0^+ \in dt, X_{\tau_0^+} \in du, -X_{\tau_0^+-} \in dv, \Delta X_{\tau_0^+} = \Delta P_{\tau_0^+}^k)$$
(3.10)
$$= \Pi_{P^k}(du + v) \widehat{\mathcal{U}}(dt, dv) \mathcal{U}(\{0\}, \{0\}).$$

REMARK 3.3 (First passage by creeping; cf. [1, 7] for details). In Theorem 3.1 and Corollary 3.2, we investigated the first passage of $X$ caused by a jump. However, $X$ may also creep over the barrier $x \in \mathbb{R}$, in which case $\mathbb{P}(X_{\tau_x^+} = x) > 0$ holds. According to [1], Theorem VI.19, is equivalent to
$$D_H = \lim_{\beta \uparrow \infty} \frac{\kappa(0, \beta)}{\beta} > 0.$$

If $X$ is of bounded variation, then $X$ creeps upwards if and only if $D_X > 0$, see [7], Section 6.4, and [19], Theorem 7.11. The linear drift $D_X$ is deterministic and so dependence between the jumps does not affect the creeping of $X$. If $X$ has a Gaussian component, then from $A = 2D_H D_{\widehat{H}}$ (see [7], Corollary 4.4(i)) $D_H > 0$ follows. So dependence between the jumps does not affect that $X$ can creep. If $X$ is of unbounded variation, but has no Gaussian component, then $X$ creeps upwards if and only if

$$\int_0^1 \frac{x\Pi_X([x,\infty))}{\int_{-x}^0 (\int_{-1}^u \Pi_X((-\infty,z])\,dz)\,du}\,dx$$
$$= \int_0^1 x(\Pi_{P^1} + \Pi_{P^2} + \Pi_{P^3} + \Pi_{P^4} + \Pi_{P^5})([x,\infty))$$
$$\times \left(\int_{-x}^0 \left(-\int_{-1}^u \overline{\Pi}_{S^{1,-}}(z) + \overline{\Pi}_{S^{2,-}}(z) + \overline{\Pi}_{P^4}^-(z)\right.\right.$$
$$\left.\left. + \overline{\Pi}_{P^5}^-(z) + \overline{\Pi}_{S^{--}}(z)\,dz\right)du\right)^{-1}dx$$

is finite. So only in this case, the dependence between the jumps can influence the possibility of creeping.

REMARK 3.4 (Quintuple law for CPPes). We work with the *weakly ascending* ladder process $(L^{-1}, H)$ and the *strictly* descending ladder process $(\widehat{L}^{-1*}, \widehat{H}^*)$ as in Section 2. We consider the *last time of the previous maximum of $X$* before time $t$ defined by $\overline{G}$ in (1.4) and the *first time of the previous minimum of $X$* before time $t$; i.e.,

(3.11) $$\underline{G}_t^* := \inf\{s < t : X_s = \underline{X}_t\},$$

see [8], Theorem 4, Remark 5, page 98, and [19], pages 167–168 and page 194. The quintuple law of Theorem 3.1 holds also for CPPes with $\widehat{\mathcal{U}}$ replaced by $\widehat{\mathcal{U}}^*$ to indicate that this is the potential measure of the strictly descending ladder process as in Section 2. The result of Corollary 3.2 changes, since $\overline{G}_{\tau_0^+-} > 0$ a.s., and we obtain for $u > 0, t \geq 0, s > 0, v \geq 0$,

$$\mathbb{P}(\tau_0^+ - \overline{G}_{\tau_0^+-} \in dt, \overline{G}_{\tau_0^+-} \in ds, X_{\tau_0^+} \in du, -X_{\tau_0^+-} \in dv, \Delta X_{\tau_0^+} = \Delta P_{\tau_0^+}^k)$$
$$= \Pi_{P^k}(du+v)\widehat{\mathcal{U}}^*(dt,dv)\mathcal{U}(ds,\{0\}), \qquad k = 1,\ldots,5.$$

The proof of the quintuple law for CPPes is analogous to Case 1 of the proof of Theorem 3.1 in Section 7. The only subtlety is in the Wiener–Hopf factorization, where we have to assign the mass given by the probabilities $\mathbb{P}(X_t = 0)$ for $t \geq 0$ to one or the other of the integrals, which define $\kappa$ and $\widehat{\kappa}$; see equations (6.19) and (6.20) in [19]. With the definition of $\underline{G}^*$ in (3.11) we assign the mass to $\kappa$; cf. [19], pages 167–168.

**4. Two explicit situations.** Whereas it is comparably easy to understand the influence of the last jump of the Lévy process, since it is independent of the past, to trace the influence of the dependence within the potential measures $\mathcal{U}$ and $\widehat{\mathcal{U}}$ given in (3.6) and (3.7) is rather involved. The ladder processes depend on the chosen local times at the maxima and minima, respectively, which in general can not be written as functionals of the path of $X$. In this section, we present two situations of Theorem 3.1, where 0 is irregular for $(0,\infty)$ and $X$ is spectrally positive.

4.1. *Spectrally positive CPP.* Let $(S^1, S^2)$ be a spectrally positive CPP and $X$ be given by

$$X_t = S^1_t + S^2_t, \qquad t \geq 0, \tag{4.1}$$

and let $\lambda > 0$ denote the jump intensity of $X$ and $F$ the d.f. of the i.i.d. jump sizes of $X$; note that both are determined by the marginal frequencies, marginal jump sizes and the dependence structure; cf. [10] for details. Set $W_0 = 0$ and denote by $(W_n)_{n \in \mathbb{N}}$ the arrival times of the jumps of $X$. $(W_n)_{n \in \mathbb{N}_0}$ is a renewal process, whose i.i.d. increments are expo($\lambda$)-distributed. Then $(X_{W_n})_{n \in \mathbb{N}_0}$ is a random walk and we can apply the result of Theorem 2.1. Because of $X = \overline{X}$ a.s. and $\overline{G}_t = \sup\{s < t : \overline{X}_s = X_s\} = t$ for all $t \geq 0$, we have $x - \overline{X}_{\tau_x^+ -} = x - X_{\tau_x^+ -}$ and $\overline{G}_{\tau_x^+ -} = \tau_x^+$ a.s. and the quintuple law reduces to a triple law. Recall the definition of the convolution $F^{n*}$ before Theorem 2.5.

THEOREM 4.1. *Suppose $X$ is given by (4.1). Consider first passage of $X$ over a constant barrier $x > 0$. Then for $u > 0, v \in [0,x], s > 0$ and for $k = 1, 2, 3$,*

$$\mathbb{P}(\tau_x^+ \in ds, X_{\tau_x^+} - x \in du, x - X_{\tau_x^+ -} \in dv, \Delta X_{\tau_x^+} = \Delta P^k_{\tau_x^+})$$

$$= \Pi_{P^k}(du + v) \sum_{n=0}^{\infty} \frac{(\lambda s)^n}{n!} e^{-\lambda s}\, ds\, F^{n*}(x - dv).$$

*By construction of the $P^k$ as being independent, the d.f. $F$ has representation*

$$F = \frac{1}{\lambda} \sum_{k=1}^{3} \lambda_{P^k} F_{P^k}, \tag{4.2}$$

*where the $\lambda_{P^k}$ denote the jump intensities of $P^k$ defined in (3.4).*

For the barrier $x = 0$, the result reduces even further.

COROLLARY 4.2. *Suppose $X$ is given as in (4.1). Consider the first passage of $X$ over the barrier $x = 0$. Then $-X_{\tau_0^+ -} = 0$ a.s. and for $u > 0$, $s > 0$ and for $k = 1, 2, 3$,*

$$\mathbb{P}(\tau_0^+ \in ds, X_{\tau_0^+} \in du, \Delta X_{\tau_0^+} = \Delta P^k_{\tau_0^+}) = \Pi_{P^k}(du) e^{-\lambda s}\, ds.$$

4.2. *Subordinator with negative drift and finite mean.* Let $(S^1, S^2)$ be a (driftless) subordinator and $X$ be given by

(4.3) $$X_t = S_t - ct = S_t^1 + S_t^2 - ct, \qquad t \geq 0,$$

with negative drift $D_X = -c < 0$. We denote its Lévy measure by $\Pi_S$ and recall the characteristic exponent of $X$ from (1.1), which reduces to

(4.4) $$\Psi_X(\theta) = \Psi_S(\theta) + ic\theta = \int_0^\infty (1 - e^{i\theta x}) \Pi_S(dx) + ic\theta, \qquad \theta \in \mathbb{R}.$$

Further, we suppose

(4.5) $$0 < \mathbb{E}[S_1] = \mu_S = \int_0^\infty x \Pi_S(dx) < c < \infty,$$

such that $\lim_{t \to \infty} X_t = -\infty$ a.s.

Under these conditions, the ascending ladder process $(L^{-1}, H)$ of $X$ is a killed bivariate CPP, and we denote its jump size distribution by $F_{\mathcal{L}^{-1}\mathcal{H}}(ds, dx)$. We denote by $F_{\mathcal{L}^{-1}\mathcal{H}}^{n*}$ the $n$-fold bivariate convolution of $F_{\mathcal{L}^{-1}\mathcal{H}}$, where $F_{\mathcal{L}^{-1}\mathcal{H}}^{0*}(ds, dz) = \delta_{(0,0)}(ds, dz)$ is the Dirac-measure in $(0,0)$ and $F_{\mathcal{L}^{-1}\mathcal{H}}^{1*} = F_{\mathcal{L}^{-1}\mathcal{H}}$. Since $L^{-1}$ and $H$ always jump together, the convolution is taken componentwise; i.e., $F_{\mathcal{L}^{-1}\mathcal{H}}^{n*}$ is the distribution of the sum of $n$ independent jumps with bivariate d.f. $F_{\mathcal{L}^{-1}\mathcal{H}}$. In this situation, we can calculate $\mathcal{U}$ and $\widehat{\mathcal{U}}$ explicitly, which we do in the proof in Section 7.

THEOREM 4.3. *Suppose $X$ is given by (4.3), and that (4.5) holds. Consider first passage of $X$ over a fixed barrier $x > 0$. Then for $u > 0$, $y \in [0, x]$, $v \geq y$, $s \geq 0$, $t \geq 0$ and for $k = 1, 2, 3$,*

(4.6) $$\begin{aligned}\mathbb{P}(\tau_x^+ - \overline{G}_{\tau_x^+-} &\in dt, \overline{G}_{\tau_x^+-} \in ds, X_{\tau_x^+} - x \in du, \\ x - X_{\tau_x^+-} &\in dv, x - \overline{X}_{\tau_x^+-} \in dy, \Delta X_{\tau_x^+} = \Delta P_{\tau_x^+}^k) \\ &= \Pi_{P^k}(du + v)\mathbb{P}(\tau_{-(v-y)}^- \in dt)\,dv\,\frac{1}{c}\sum_{n=0}^\infty \left(\frac{\mu_S}{c}\right)^n F_{\mathcal{L}^{-1}\mathcal{H}}^{n*}(ds, x - dy),\end{aligned}$$

*where the bivariate jump size d.f. $F_{\mathcal{L}^{-1}\mathcal{H}}$ is given by*

(4.7) $$F_{\mathcal{L}^{-1}\mathcal{H}}(ds, dz) = \frac{1}{\mu_S} \int_{[0, \infty)} \Pi_S(dz + \theta) \mathbb{P}(\tau_{-\theta}^- \in ds)\,d\theta.$$

If we are only interested in the space variables, we can integrate out time in the above quintuple law and obtain as proved in Section 7 the following.

COROLLARY 4.4. *In the situation of Theorem 4.3, for $u > 0$, $y \in [0, x]$, $v \geq y$, we get for $k = 1, 2, 3$,*

$$\mathbb{P}(X_{\tau_x^+} - x \in du, x - X_{\tau_x^+-} \in dv, x - \overline{X}_{\tau_x^+-} \in dy, \Delta X_{\tau_x^+} = \Delta P_{\tau_x^+}^k)$$

(4.8)
$$= \Pi_{P^k}(du+v)\,dv\,\frac{1}{c}\sum_{n=0}^{\infty}\left(\frac{\mu_S}{c}\right)^n F_{\mathcal{H}}^{n*}(x-dy),$$

(4.9) $\quad \mathbb{P}(\Delta X_{\tau_x^+} = \Delta P^k_{\tau_x^+}) = \frac{1}{c}\int_0^x\int_y^{\infty}\overline{\Pi}_{P^k}(v)\,dv\sum_{n=0}^{\infty}\left(\frac{\mu_S}{c}\right)^n F_{\mathcal{H}}^{n*}(x-dy),$

(4.10) $\quad\quad \mathbb{P}(\tau_x^+ < \infty) = \left(1-\frac{\mu_S}{c}\right)\sum_{n=1}^{\infty}\left(\frac{\mu_S}{c}\right)^n \overline{F_{\mathcal{H}}^{n*}}(x).$

*The d.f. $F_{\mathcal{H}}$ is for $z > 0$ defined as*

(4.11) $\quad\quad F_{\mathcal{H}}(dz) = \frac{1}{\mu_S}\overline{\Pi}_S(z)\,dz = \frac{1}{\mu_S}(\overline{\Pi}_{P^1}+\overline{\Pi}_{P^2}+\overline{\Pi}_{P^3})(z)\,dz.$

*Here, $F_{\mathcal{H}}^{0*}(dz) = \delta_0(dz)$ and for $n \in \mathbb{N}$*

$$F_{\mathcal{H}}^{n*}(dz) = \frac{1}{\mu_S^n}\overline{\Pi}_S^{n\otimes}(z)\,dz$$

*with $\overline{\Pi}_S^{1\otimes} = \overline{\Pi}_S$ and $\overline{\Pi}_S^{2\otimes}(z) := \int_0^z \overline{\Pi}_S(z-y)\overline{\Pi}_S(y)\,dy$.*

REMARK 4.5. When the jump part $S$ in (4.3) is a CPP with jump size d.f. $F$ and $\mathbb{E}[\Delta S] = \mu_{\Delta S}$ then, under the conditions of Theorem 4.3, for $x > 0$

$$F_{\mathcal{H}}(x) = \frac{1}{\mu_{\Delta S}}\int_0^x \overline{F}(z)\,dz$$

*and (4.10) is the celebrated Pollaczek–Khintchine formula.*

COROLLARY 4.6. *Suppose $X$ is given by (4.3) and (4.5) holds. Consider the first passage of $X$ over the barrier $x = 0$. Then for $u > 0$, $v \geq 0$ and $t > 0$ and for $k = 1,2,3$,*

(4.12)
$$\mathbb{P}(\tau_0^+ \in dt, X_{\tau_0^+} \in du, -X_{\tau_0^+-} \in dv, \Delta X_{\tau_0^+} = \Delta P^k_{\tau_0^+})$$
$$= \Pi_{P^k}(du+v)\mathbb{P}(\tau_{-v}^- \in dt)\frac{1}{c}\,dv.$$

*Further,*

(4.13) $\quad \mathbb{P}(X_{\tau_0^+} \in du, -X_{\tau_0^+-} \in dv, \Delta X_{\tau_0^+} = \Delta P^k_{\tau_0^+}) = \Pi_{P^k}(du+v)\frac{1}{c}dv,$

(4.14) $\quad\quad\quad\quad\quad\quad \mathbb{P}(\Delta X_{\tau_0^+} = \Delta P^k_{\tau_0^+}) = \frac{\mu_{P^k}}{c},$

(4.15) $\quad\quad\quad\quad\quad\quad\quad\quad \mathbb{P}(\tau_0^+ < \infty) = \frac{\mu_S}{c}.$

The identities (4.13) and (4.14) are generalizations of Theorem 2.2(i) in [13], where only independence is treated. When we compare (4.10) and (4.15), we see that the ruin probability for the barrier $x=0$ is not affected by the dependence in contrast to barriers $x>0$.

**5. Dependence modeling by a Lévy copula.** Our goal is to study the effect of dependence between the jumps of $X^1$ and $X^2$ for the quintuple law. This dependence affects the bivariate potential measures $\mathcal{U}, \widehat{\mathcal{U}}$ and also the factors $\Pi_{\mathcal{P}^k}$. As dependence structure, we use the concept of a Lévy copula.

5.1. *Lévy copulas.* Similarly to copulas for d.f.s, Lévy copulas have been suggested to model the dependence in the jump behaviour of Lévy processes. The basic idea is to invoke an analogue of Sklar's theorem (2.4) for Lévy measures (cf. [15], Theorem 3.6) and model the dependence structure by a Lévy copula on the marginal Lévy measures. Since, with the exception of a CPP, all Lévy measures have a singularity in 0, we follow Kallsen and Tankov [15] and introduce a Lévy copula on the *tail integral*, which is defined for each quadrant separately. For convenience, we set $\overline{\mathbb{R}} := (-\infty, \infty]$.

DEFINITION 5.1 (Bivariate tail integral and its margins).

(1) Let $(X^1, X^2)$ be an $\mathbb{R}^2$-valued Lévy process with Lévy measure $\Pi$. Its *tail integral* is the function $\overline{\Pi} : (\mathbb{R} \setminus \{0\})^2 \to \mathbb{R}$ defined quadrantwise as

$$\overline{\Pi}(x_1, x_2) := \begin{cases} \Pi((x_1, \infty) \times (x_2, \infty)), & x_1 > 0,\ x_2 > 0, \\ -\Pi((x_1, \infty) \times (-\infty, x_2]), & x_1 > 0,\ x_2 < 0, \\ -\Pi((-\infty, x_1] \times (x_2, \infty)), & x_1 < 0,\ x_2 > 0, \\ \Pi((-\infty, x_1] \times (-\infty, x_2]), & x_1 < 0,\ x_2 < 0. \end{cases}$$

(2) The *marginal tail integrals* $\overline{\Pi}_1$ and $\overline{\Pi}_2$ are the tail integrals of the processes $X^1$ and $X^2$, respectively, and we set for $i=1,2$

$$\overline{\Pi}_i(x) := \begin{cases} \Pi_i((x, \infty)), & x > 0, \\ -\Pi_i((-\infty, x]), & x < 0. \end{cases}$$

THEOREM 5.2 (Sklar's theorem for Lévy copulas, [15], Theorem 3.6).

(1) *Let $(X^1, X^2)$ be an $\mathbb{R}^2$-valued Lévy process. We can write the joint tail integral of $(X^1, X^2)$ as*

$$(5.1) \qquad \overline{\Pi}(x_1, x_2) = \widehat{C}(\overline{\Pi}_1(x_1), \overline{\Pi}_2(x_2)), \qquad (x_1, x_2) \in (\mathbb{R} \setminus \{0\})^2,$$

*where the Lévy copula $\widehat{C}$ is unique on $\overline{\mathrm{Ran}\,\overline{\Pi}_1} \times \overline{\mathrm{Ran}\,\overline{\Pi}_2}$.*

(2) *Let $\widehat{C}$ be a bivariate Lévy copula and $\overline{\Pi}_1$, $\overline{\Pi}_2$ tail integrals of Lévy processes. Then there exists an $\mathbb{R}^2$-valued Lévy process $(X^1, X^2)$ whose components have tail integrals $\overline{\Pi}_1$, $\overline{\Pi}_2$ satisfying (5.1) for every $(x_1, x_2) \in (\mathbb{R} \setminus \{0\})^2$. The Lévy measure $\Pi$ of $(X^1, X^2)$ is uniquely determined by $\widehat{C}$ and $\overline{\Pi}_1$, $\overline{\Pi}_2$.*

Various possibilities to construct parametric Lévy copula families can be found in the literature. Proposition 5.5 of [6] shows how to construct a positive Lévy copula from a distributional copula, Theorem 5.1 of [15] presents Lévy copulas as limits of distributional copulas, and their Theorem 6.1 shows how to construct Archimedean Lévy copulas analogously to Archimedean copulas. We present some examples for later use.

EXAMPLE 5.3. (a) The *independence Lévy copula* is defined as (cf. [15], Proposition 4.1)

$$(5.2) \qquad \widehat{C}_\perp(u,v) = u 1_{\{v=\infty\}} + v 1_{\{u=\infty\}}.$$

(b) The *complete dependence Lévy copula* is defined as (cf. [15], Theorem 4.4)

$$(5.3) \qquad \widehat{C}_\|(u,v) = \min\{|u|,|v|\} 1_K(u,v) \operatorname{sgn}(u) \operatorname{sgn}(v),$$

where $K := \{x \in \mathbb{R}^2 : \operatorname{sgn}(x_1) = \operatorname{sgn}(x_2)\} = (-\infty, 0)^2 \cup [0, \infty)^2$.

(c) For $\theta > 0$ and $\eta \in [0,1]$, the *Clayton–Lévy copula* is defined as [cf. [15], equation (6.5)]

$$(5.4) \qquad \widehat{C}_{\eta,\theta}(u,v) = (|u|^{-\theta} + |v|^{-\theta})^{-1/\theta} (\eta 1_{\{uv \geq 0\}} - (1-\eta) 1_{\{uv < 0\}}).$$

For $\eta = 1$, the two components always jump in the same direction, for $\eta = 0$ in opposite direction. The parameter $\theta$ models the degree of dependence: for $\eta = 1$ and $\theta \to 0$, we obtain the independence model and for $\theta \to \infty$ the complete dependence model. Further, Clayton–Lévy copulas are homogeneous of order 1, i.e., $\widehat{C}_{\eta,\theta}(cu, cv) = c \widehat{C}_{\eta,\theta}(u,v)$ for all $c > 0$, and continuous on $\overline{\mathbb{R}}^2$. This Lévy copula has been used frequently, see e.g., Example 5.9 in [6] and [4, 5, 10].

(d) For $\zeta > 0$ and $\eta \in (0,1)$, a nonhomogeneous Archimedean Lévy copula can be found by invoking Theorem 6.1 of [15]. Defining $\varphi = \varphi_\zeta : [-1,1] \to [-\infty, \infty]$ by

$$\varphi_\zeta(x) = \eta \zeta \frac{x}{1-x} 1_{\{x \geq 0\}} - (1-\eta) \zeta \frac{|x|}{1-|x|} 1_{\{x < 0\}}$$

yields the nonhomogeneous Lévy copula

$$(5.5) \qquad \widehat{C}_{\eta,\zeta}(u,v) = \frac{|uv|}{|u|+|v|+\zeta}(\eta 1_{\{uv \geq 0\}} - (1-\eta) 1_{\{uv < 0\}}),$$

which is continuous on $\overline{\mathbb{R}}^2$. From $\lim_{\zeta \to \infty} \widehat{C}_{\eta,\zeta}(u,v) = 0$ for all $u, v \in \mathbb{R}$ and for all $\eta \in (0,1)$, we obtain the independence Lévy copula, and for $\zeta \to 0$, we obtain the Clayton–Lévy copula with parameter $\theta = 1$.

5.2. *Calculating the quantities in the quintuple law.* Analogously to Definition 5.1(2), we define the tail integrals $\overline{\Pi}_{P^k}$ for $k = 1, \ldots, 5$ of the single and joint jump components. The following result shows the influence of a specific Lévy copula on the tail integrals, where we set $\overline{\Pi}_{S^2}(0) := \Pi_{S^2}((0, \infty))$. A proof can be found in Section 7.

THEOREM 5.4. *Suppose that the jump parts $S^1$ and $S^2$, given in (3.1), have absolutely continuous Lévy measures $\Pi_{S^i}$ and the dependence between their jumps is modeled by a twice continuously differentiable Lévy copula $\widehat{C}$. Then the tail integrals in Theorem 3.1 and Corollary 3.2, respectively, are given for $z > 0$,*

$$(5.6) \quad \overline{\Pi}_{P^1}(z) = \overline{\Pi}_{S^1}(z) - \lim_{y \downarrow 0} \widehat{C}(\overline{\Pi}_{S^1}(z), \overline{\Pi}_{S^2}(y)) + \lim_{y \uparrow 0} \widehat{C}(\overline{\Pi}_{S^1}(z), \overline{\Pi}_{S^2}(y)),$$

$$(5.7) \quad \overline{\Pi}_{P^2}(z) = \overline{\Pi}_{S^2}(z) - \lim_{x \downarrow 0} \widehat{C}(\overline{\Pi}_{S^1}(x), \overline{\Pi}_{S^2}(z)) + \lim_{x \uparrow 0} \widehat{C}(\overline{\Pi}_{S^1}(x), \overline{\Pi}_{S^2}(z)),$$

$$(5.8) \quad \overline{\Pi}_{P^3}(z) = \int_0^\infty \left. \frac{\partial \widehat{C}(u,v)}{\partial u} \right|_{u=\overline{\Pi}_{S^1}(x), v=\overline{\Pi}_{S^2}((z-x)\vee 0)} \Pi_{S^1}(dx),$$

$$\overline{\Pi}_{P^4}(z) = \int_z^\infty \left[ \left. \frac{\partial \widehat{C}(u,v)}{\partial u} \right|_{u=\overline{\Pi}_{S^1}(x)} \right]_{\lim_{a \uparrow 0} \overline{\Pi}_{S^2}(a)}^{\overline{\Pi}_{S^2}(z-x)} \Pi_{S^1}(dx),$$

$$\overline{\Pi}_{P^5}(z) = \int_z^\infty \left[ \left. \frac{\partial \widehat{C}(u,v)}{\partial v} \right|_{v=\overline{\Pi}_{S^2}(y)} \right]_{\lim_{a \uparrow 0} \overline{\Pi}_{S^1}(a)}^{\overline{\Pi}_{S^1}(z-y)} \Pi_{S^2}(dy).$$

*If $\widehat{C}$ is left-continuous in the second coordinate in $\infty$ and $\Pi_{S^2}((0, \infty)) = \Pi_{S^2}((-\infty, 0)) = \infty$, then (5.6) reduces to $\Pi_{P^1} \equiv 0$.*

*The analogous result holds for the first coordinate with $\Pi_{P^2} \equiv 0$.*

The following result is a simple consequence of Theorem 5.4.

COROLLARY 5.5 ([5], Proposition 2.16). *Assume that the conditions of Theorem 5.4 hold and that $S^1$ and $S^2$ are spectrally positive. Then the tails (5.6)–(5.8) reduce to*

$$\overline{\Pi}_{P^1}(z) = \overline{\Pi}_{S^1}(z) - \lim_{y \downarrow 0} \widehat{C}(\overline{\Pi}_{S^1}(z), \overline{\Pi}_{S^2}(y)),$$

$$\overline{\Pi}_{P^2}(z) = \overline{\Pi}_{S^2}(z) - \lim_{x \downarrow 0} \widehat{C}(\overline{\Pi}_{S^1}(x), \overline{\Pi}_{S^2}(z)),$$

$$\overline{\Pi}_{P^3}(z) = \int_0^\infty \left. \frac{\partial \widehat{C}(u,v)}{\partial u} \right|_{u=\overline{\Pi}_{S^1}(x), v=\overline{\Pi}_{S^2}(0 \vee (z-x))} \Pi_{S^1}(dx).$$

If $\widehat{C}$ is left-continuous in the second coordinate in $\infty$ and $\Pi_{S^2}((0,\infty)) = \infty$, then $\Pi_{P^1} \equiv 0$.

The analogous result holds for the first coordinate with $\Pi_{P^2} \equiv 0$.

The following lemma shows that for a left-continuous and homogeneous Lévy copula $\widehat{C}$ single jumps always have a marginal lighter tail than the corresponding component.

LEMMA 5.6. *Assume that the conditions of Corollary 5.5 hold. If $\widehat{C}$ is left-continuous in the $j$th coordinate in $\infty$ and homogeneous, then for $i \neq j$,*

$$\overline{\Pi}_{P^i}(z) = o(\overline{\Pi}_{S^i}(z)), \qquad z \to \infty. \tag{5.9}$$

Now, we apply our results to the situations of Section 4.

THEOREM 5.7. *Suppose that the jump parts $S^1$ and $S^2$, given in (3.1), have absolutely continuous Lévy measures $\Pi_{S^i}$ and the dependence between their jumps is modeled by a twice continuously differentiable Lévy copula $\widehat{C}$.*

*(1) In the situation of Section 4.1, when $S^1$ and $S^2$ are spectrally positive CPPes with jump intensities $\lambda_1$, $\lambda_2$ and jump size d.f.s $F_1$, $F_2$, Theorem 4.1 and Corollary 4.2 hold with*

$$\lambda = \lambda_1 + \lambda_2 - \widehat{C}(\lambda_1, \lambda_2), \tag{5.10}$$

*and for $z > 0$*

$$\overline{F}(z) = \frac{1}{\lambda}\bigg(\lambda_1 \overline{F}_1(z) - \widehat{C}(\lambda_1 \overline{F}_1(z), \lambda_2) + \lambda_2 \overline{F}_2(z) - \widehat{C}(\lambda_1, \lambda_2 \overline{F}_2(z))$$
$$+ \lambda_1 \int_0^\infty \frac{\partial \widehat{C}(u,v)}{\partial u}\bigg|_{u=\lambda_1 \overline{F}_1(x), v=\lambda_2 \overline{F}_2(0 \vee (z-x))} F_1(dx)\bigg). \tag{5.11}$$

*(2) In the situation of Section 4.2, when $S = S^1 + S^2$ is a subordinator, Corollary 4.4 holds for $F_{\mathcal{H}}(z)$ given for $z > 0$ by*

$$F_{\mathcal{H}}(dz) = \frac{1}{\mu_S}\bigg(\overline{\Pi}_{S^1}(z) - \lim_{y \downarrow 0} \widehat{C}(\overline{\Pi}_{S^1}(z), \overline{\Pi}_{S^2}(y))$$
$$+ \overline{\Pi}_{S^2}(z) - \lim_{x \downarrow 0} \widehat{C}(\overline{\Pi}_{S^1}(x), \overline{\Pi}_{S^2}(z))$$
$$+ \int_0^\infty \frac{\partial \widehat{C}(u,v)}{\partial u}\bigg|_{u=\overline{\Pi}_{S^1}(x), v=\overline{\Pi}_{S^2}(0 \vee (z-x))} \Pi_{S^1}(dx)\bigg) dz.$$

REMARK 5.8 (Comparison of random walk and Lévy process modeling). Let $(X^1, X^2)$ be a spectrally positive CPP (without drift) with marginal intensities $\lambda_1, \lambda_2$ and absolutely continuous marginal jump size d.f.s $F_1, F_2$. Denote by $(W_n^i)_{n \in \mathbb{N}_0}$ the arrival times of the jumps of $X^i$. We use the embedded random walk structure by defining $Z_n^i := X_{W_n^i}^i$. Then the d.f. of the increments of $Z^i$ is equal to $F_i$ and, by absolute continuity, $Z^1$ and $Z^2$ always jump together. If we model the dependence between jumps by a distributional copula $C$, then with equation (2.5) the tail of the jump size d.f. $F_Z$ of $Z = Z^1 + Z^2$ is given by

$$\overline{F}_Z(z) = \int_0^\infty \left(1 - \frac{\partial C(u,v)}{\partial u}\bigg|_{u=F_1(x), v=F_2((z-x)\vee 0)}\right) F_1(dx).$$

Rewriting this expression in terms of the distributional survival copula $\widetilde{C}(u,v) := u + v - 1 + C(1-u, 1-v)$ (see [20], equation (2.6.2)) yields

$$(5.12) \qquad \overline{F}_Z(z) = \int_0^\infty \frac{\partial \widetilde{C}(u,v)}{\partial u}\bigg|_{u=\overline{F}_1(x), v=\overline{F}_2((z-x)\vee 0)} F_1(dx).$$

When we consider, however, the Lévy process $(X^1, X^2)$ and use a Lévy copula $\widehat{C}$, then the tail $\mathbb{P}(\Delta X^1 + \Delta X^2 > z)$ is given by (5.11). Comparing (5.12) and (5.11), the most apparent differences are the first four summands in (5.11). These summands represent the possibility of single jumps of $X^i$. They are missing in (5.12) since the random walks $Z^1$ and $Z^2$ always jump together by construction. But also the last integrals in (5.11) and (5.12), which represent the common jumps of $X^1$ and $X^2$, differ. Furthermore, we have obviously different spaces: a distributional survival copula $\widetilde{C} \colon [0,1]^2 \to [0,1]$ and a Lévy copula $\widehat{C} \colon (-\infty, \infty]^2 \to (-\infty, \infty]$.

5.3. *Examples for different dependence structures.* We present four examples for different dependence structures, modeled by a Lévy copula $\widehat{C}$, and characterize all quantities of Theorem 4.1 and Corollary 4.4.

5.3.1. *Independence.* Suppose $S^1$ and $S^2$ are independent, i.e., $S^1$ and $S^2$ a.s. never jump together, cf. Example 12.10(i) of Sato [21]. Then $P^1 = S^{1,+} = S^1$, $P^2 = S^{2,+} = S^2$ and $P^3 = S^{1,++} + S^{2,++} = 0$. This corresponds to Example 5.3(a), and we obtain

$$\Pi_S(dz) = (\Pi_{S^1} + \Pi_{S^2})(dz) = (\Pi_{P^1} + \Pi_{P^2})(dz).$$

In the situation of Section 4.1, when the jump parts $S^1$ and $S^2$ are spectrally positive CPPes with intensities $\lambda_1$ and $\lambda_2$ and jump size d.f.s $F_1$ and $F_2$, we get in Theorem 4.1 and in Corollary 4.2 the identities $\lambda = \lambda_1 + \lambda_2$ and

$$F(dz) = \frac{1}{\lambda}\Pi_S(dz) = \left(\frac{\lambda_1}{\lambda_1 + \lambda_2}F_1 + \frac{\lambda_2}{\lambda_1 + \lambda_2}F_2\right)(dz).$$

In the situation of Section 4.2, when $X$ is a subordinator with negative drift, Corollary 4.4 holds with

$$F_{\mathcal{H}}(dz) = \frac{1}{\mu_S}(\overline{\Pi}_{S^1}(z) + \overline{\Pi}_{S^2}(z))\, dz.$$

In this case (4.13) and (4.14) are the results of [13], Theorem 2.2(i).

5.3.2. *Complete dependence.* The jumps of $(S^1, S^2)$ are said to be completely dependent, if there is a strictly ordered subset $\mathcal{S} \subset [0, \infty)^2$ such that $(\Delta S_t^1, \Delta S_t^2) \in \mathcal{S}$ for every $t > 0$ (except for some null set of paths), see [15], Definition 4.2. This means that $S^1$ and $S^2$ a.s. jump together and so $P^1 = P^2 = 0$ and $P^3 = S^1 + S^2 = S$. This is the situation of Example 5.3(b). In Section 4.1, we get in Theorem 4.1 and Corollary 4.2 the identities $\lambda = \lambda_1 = \lambda_2$ and

$$F(dz) = \mathbb{P}(\Delta S^1 + \Delta S^2 \in dz),$$

where the jump sizes $\Delta S^1$ and $\Delta S^2$ are completely dependent. In Corollary 4.4, we get $F_{\mathcal{H}}(dz) = \overline{\Pi}_{P^3}(z)\, dz/\mu_S$. A simple example for complete dependence is $S^1 \equiv S^2$, then $F(dz) = F_1(dz/2)$ and

$$F_{\mathcal{H}}(dz) = \frac{1}{\mu_S}\overline{\Pi}_{S^1}(z/2)\, dz.$$

5.3.3. *Clayton–Lévy copula.* In both situations of Section 4, the jump part $S$ has only positive jumps and so we must have $\eta = 1$ in (5.4). By absolute continuity, the tail integral of $S$ is for $z > 0$ given by

$$\begin{aligned}
\overline{\Pi}_S(z) &= \overline{\Pi}_{P^1}(z) + \overline{\Pi}_{P^2}(z) + \overline{\Pi}_{P^3}(z) \\
&= \overline{\Pi}_{S^1}(z) - (\overline{\Pi}_{S^1}(z)^{-\theta} + (\Pi_{S^2}((0,\infty)))^{-\theta})^{-1/\theta} \\
&\quad + \overline{\Pi}_{S^2}(z) - (\overline{\Pi}_{S^2}(z)^{-\theta} + (\Pi_{S^1}((0,\infty)))^{-\theta})^{-1/\theta} \\
&\quad + \int_{(0,\infty)} (\overline{\Pi}_{S^1}(x)^{-\theta} + \overline{\Pi}_{S^2}((z-x) \vee 0)^{-\theta})^{-1/\theta - 1} \\
&\quad \times \overline{\Pi}_{S^1}(x)^{-\theta - 1} \Pi_{S^1}(dx).
\end{aligned} \quad (5.13)$$

If $\Pi_{S^1}((0,\infty)) = \infty$, then we see from (5.13) that $\Pi_{P^2} = 0$ holds, i.e., $S^2$ has no single jumps. So if $\Pi_{S^1}$ and $\Pi_{S^2}$ are infinite measures, then there are infinitely many common jumps and no single jumps. If $\Pi_{S^1}((0,\infty)) = \infty$ and $\Pi_{S^2}((0,\infty)) = \lambda_2 < \infty$, then the intensity rate of the common jumps reduce to $\Pi((0,\infty) \times (0,\infty)) = \lim_{a \to \infty} \widehat{C}_\theta(a, \lambda_2) = \lambda_2$. If $(S^1, S^2)$ is a CPP, then we get the result of [5], Proposition 3.1. In Section 4.1, Theorem 4.1 holds with

$$\lambda = \lambda_1 + \lambda_2 - (\lambda_1^{-\theta} + \lambda_2^{-\theta})^{-1/\theta} \quad \text{and} \quad \overline{F}(z) = \frac{1}{\lambda}\overline{\Pi}_S(z), \qquad z > 0,$$

and in Section 4.2, Corollary 4.4 holds with $F_{\mathcal{H}}(dz) = \overline{\Pi}_S(dz)/\mu_S$. For all $u$, $v > 0$,

$$\frac{\partial \widehat{C}_\theta}{\partial \theta}(u,v) = \theta^{-2}(u^{-\theta} + v^{-\theta})^{-1/\theta - 1}(u^{-\theta}(\ln(u^{-\theta} + v^{-\theta}) + \theta \ln u)$$
$$+ v^{-\theta}(\ln(u^{-\theta} + v^{-\theta}) + \theta \ln v)) \geq 0$$

and $\overline{\Pi}_{P^1}(z) = \overline{\Pi}_{S^1}(z) - \widehat{C}_\theta(\overline{\Pi}_{S^1}(z), \Pi_{S^2}((0, \infty)))$. From this, we see that increasing the dependence parameter $\theta$ reduces the number of single jumps, see Figure 1. In the special case of two CPPes with the same marginal Lévy measures, which are exponential, i.e., $\overline{\Pi}_{S^1}(x) = \overline{\Pi}_{S^2}(x) = e^{-ax}$ for some $a > 0$, and $\theta = 1$ we find (cf. [5], Example 3.11) for $z > 0$,

$$\overline{\Pi}_S(z) = e^{-az}\left(\frac{2}{e^{az} + 1} + \frac{1}{e^{-az} + 1}\right)$$
$$+ \frac{1}{2}e^{-(1/2)az}(\arctan(e^{(1/2)az}) - \arctan(e^{-(1/2)az}))$$
$$= \frac{3 + 2e^{-az} + e^{az}}{(e^{az} + 1)(e^{-az} + 1)}$$
$$+ \frac{1}{2}e^{-(1/2)az}(\arctan(e^{(1/2)az}) - \arctan(e^{-(1/2)az}))$$
$$\sim e^{-az}\left(1 + \frac{\pi}{2}e^{-(1/2)az}\right), \qquad z \to \infty.$$

This implies that asymptotically for large $z$ the joint jumps dominate.

5.3.4. *Nonhomogeneous Archimedean Lévy copula.* In the spectrally positive situations of Section 4, we have $\eta = 1$ and for $z > 0$,

$$\overline{\Pi}_S(z) = \overline{\Pi}_{S^1}(z)\left(1 - \frac{\lambda_2}{\overline{\Pi}_{S^1}(z) + \lambda_2 + \zeta}\right)$$
(5.14)
$$+ \overline{\Pi}_{S^2}(z)\left(1 - \frac{\lambda_1}{\overline{\Pi}_{S^2}(z) + \lambda_1 + \zeta}\right)$$
$$+ \int_0^\infty \frac{\overline{\Pi}_{S^2}(0 \vee (z-x))^2 + \zeta \overline{\Pi}_{S^2}(0 \vee (z-x))}{(\overline{\Pi}_{S^1}(x) + \overline{\Pi}_{S^2}(0 \vee (z-x)) + \zeta)^2} \Pi_{S^1}(dx).$$

In Section 4.1, Theorem 4.1 and Corollary 4.2 hold with

$$\lambda = \lambda_1 + \lambda_2 - \frac{\lambda_1 \lambda_2}{\lambda_1 + \lambda_2 + \zeta} \quad \text{and} \quad \overline{F}(z) = \frac{1}{\lambda}\overline{\Pi}_S(z), \qquad z > 0,$$

and in Section 4.2, Corollary 4.4 holds with $F_{\mathcal{H}}(dz) = 1/\mu_S \overline{\Pi}_S(z)\, dz$. Contrary to the Clayton–Lévy copula, increasing the dependence parameter $\zeta$ yields that the tail integrals of the single jumps increase and those of the common jumps decrease. Further, the jump intensity $\lambda$ increases due to more single jumps. Since Lemma 5.6 does not cover the case of the non-homogeneous Lévy copula for the CPP, we calculate the following quotient explicitly:

$$\frac{\overline{\Pi}_{P^1}(z)}{\overline{\Pi}_{S^1}(z)} = 1 - \lim_{y \downarrow 0} \frac{\overline{\Pi}_{S^2}(y)}{\overline{\Pi}_{S^1}(z) + \overline{\Pi}_{S^2}(y) + \zeta}$$

$$= 1 - \frac{\lambda_2}{\overline{\Pi}_{S^1}(z) + \lambda_2 + \zeta} \to \frac{\zeta}{\lambda_2 + \zeta}, \qquad z \to \infty.$$

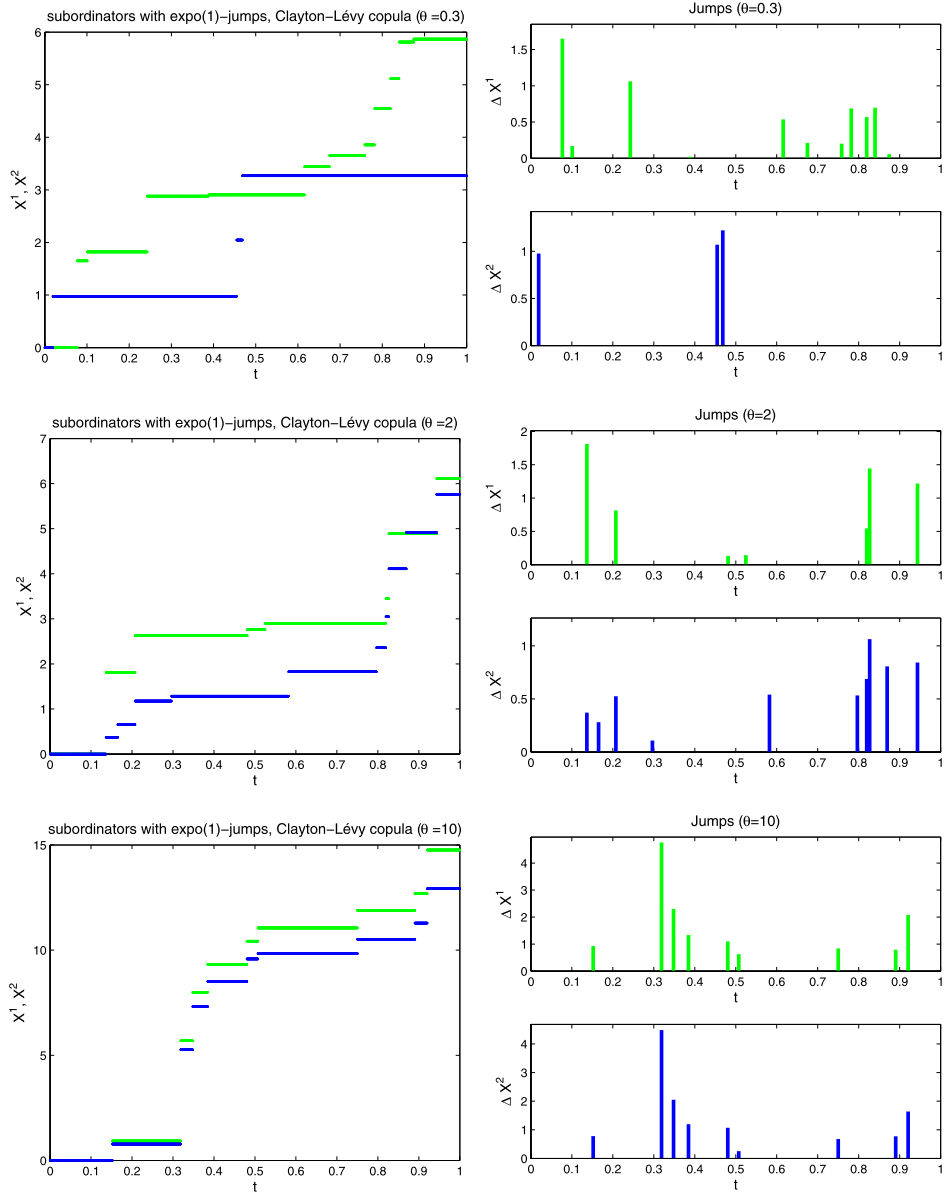

Fig. 1. *Simulations of CPPes $(X^1, X^2)$ with $\mathrm{expo}(1)$-distributed jump sizes and dependence modeled by a Clayton–Lévy copula for $\theta = 0.3$, $\theta = 2$ and $\theta = 10$; left-hand side: sample paths of CPP, right-hand side: CPP as marked point process. When $\theta$ increases, then the number of single jumps $\Delta P^1, \Delta P^2$, cf. (3.4), decreases and the number of common jumps $\Delta P^3$ increases. Further, for increasing $\theta$, the dependence between jump sizes of $X^1$ and $X^2$ increases.*

In contrast to the homogeneous Lévy copulas, the single jump Lévy measures are tail-equivalent to the Lévy measures of the components. Consequently, the Lévy measure of the sum process can or cannot be dominated by the common jumps. We present two different examples.

Take two CPPes with the same exponential marginal Lévy measures given by $\overline{\Pi}_{S^1}(x) = \overline{\Pi}_{S^2}(x) = e^{-ax}$ for some $a > 0$ and $\zeta > 2$. Then for $z > 0$ we get

$$\overline{\Pi}_{P^1}(z) = \overline{\Pi}_{P^2}(z) = e^{-az} \frac{e^{-az} + \zeta}{e^{-az} + 1 + \zeta} \sim \frac{\zeta}{1+\zeta} e^{-az}, \qquad z \to \infty.$$

It remains to calculate $\overline{\Pi}_{P^3}$:

$$\overline{\Pi}_{P^3}(z) = a \int_0^z \frac{e^{-2a(z-x)} + \zeta e^{-a(z-x)}}{(e^{-ax} + e^{-a(z-x)} + \zeta)^2} e^{-ax} \, dx$$

$$+ a \int_z^\infty \frac{1+\zeta}{(e^{-ax} + 1 + \zeta)^2} e^{-ax} \, dx =: I(z) + II(z).$$

We substitute $y = e^{ax}$ and calculate both integrals separately.

$$II(z) = \int_{e^{az}}^\infty \frac{1+\zeta}{(1+y(1+\zeta))^2} \, dy = \frac{1}{1 + e^{az}(1+\zeta)} \sim \frac{1}{1+\zeta} e^{-az}, \qquad z \to \infty.$$

For $\zeta^2 - 4e^{-az} > 0$, we get

$$I(z) = \int_1^{e^{az}} \frac{e^{-2az}y^2 + \zeta e^{-az}y}{(1 + e^{-az}y^2 + \zeta y)^2} \, dy$$

$$= \frac{e^{-az}\zeta(1 - e^{-az})}{(4e^{-az} - \zeta^2)(1 + \zeta + e^{-az})} + \frac{e^{-az}(2e^{-az} - \zeta^2)}{(4e^{-az} - \zeta^2)\sqrt{\zeta^2 - 4e^{-az}}}$$

$$\times \ln\left(\frac{(2+\zeta-\sqrt{\zeta^2-4e^{-az}})(2e^{-az}+\zeta+\sqrt{\zeta^2-4e^{-az}})}{(2+\zeta+\sqrt{\zeta^2-4e^{-az}})(2e^{-az}+\zeta-\sqrt{\zeta^2-4e^{-az}})}\right).$$

Applying l'Hospital's lemma to the last term yields $I(z) \sim (a/\zeta)ze^{-az}$ as $z \to \infty$ and, hence,

$$\overline{\Pi}_S(z) \sim \frac{a}{\zeta} z e^{-az} \qquad \text{as } z \to \infty.$$

This implies that for large $z$ the joint jumps dominate.

As a heavy-tailed example, we consider standard Pareto margins, i.e., $\overline{\Pi}_{S^1}(x) = \overline{\Pi}_{S^2}(x) = x^{-1}$ for $x \geq 1$. Then we get for $z > 1$

$$\overline{\Pi}_{P^1}(z) = \overline{\Pi}_{P^2}(z) = \frac{\zeta + z^{-1}}{1 + z(1 + \zeta)} \sim \frac{\zeta}{1 + \zeta} z^{-1} \quad \text{as } z \to \infty,$$

and for $z > 2$

$$\overline{\Pi}_{P^3}(z) = \frac{2z^2\zeta + 6z - 2z\zeta - 4}{(4 + z\zeta)(-\zeta + z\zeta + z)z}$$
$$+ \frac{2(2 + z\zeta)}{(4 + z\zeta)z\sqrt{z\zeta(4 + z\zeta)}} \ln\left(\left|\frac{z\zeta - 2\zeta + \sqrt{z\zeta(4 + z\zeta)}}{z\zeta - 2\zeta - \sqrt{z\zeta(4 + z\zeta)}}\right|\right)$$
$$\sim \frac{2}{1 + \zeta} z^{-1} \quad \text{as } z \to \infty.$$

Hence,

$$\overline{\Pi}_S(z) \sim \frac{2 + \zeta}{1 + \zeta} z^{-1} \quad \text{as } z \to \infty$$

and in contrast to the light-tailed example, common and single jumps determine $\overline{\Pi}_S(z)$ for large $z$.

**6. Applications to insurance risk theory.** We consider the situation as in [5], where the total insurance risk process of an insurance company is modelled as the sum of the components of a $d$-dimensional spectrally positive compound Poisson risk process of different business lines, which may be dependent. Our present model allows for a general spectrally positive Lévy process. This situation has been considered in [17, 18] for a one-dimensional risk model.

The company's total risk process $X$ describes the net balance of the insurance business given (for $d = 2$) by $X = X^1 + X^2$, which may involve a Gaussian baseline component and a drift, usually due to a linear premium income. The jump part

$$S = S^1 + S^2$$

of $X$ models the total amount of claims in all business lines. Ruin of the company happens, if the first hitting time $\tau_x^+$ given in (1.3) is finite. We suppose throughout this section that the net profit condition

(6.1) $$\lim_{t \to \infty} X_t = -\infty \quad \text{a.s.}$$

holds. Then the probability of first upwards passage over the barrier $x$ decreases to 0, when $x \uparrow \infty$, and our results allow us a very precise description of the ruin event caused by a jump.

Since $X$ is spectrally positive, we can choose the descending ladder process $(\widehat{L}^{-1}, \widehat{H})$ such that $\widehat{\mathcal{U}}(dx) = dx$ and, under the normalization condition (3.8), we get (as in the proof of Theorem 4.3) equation (7.3); i.e.,

$$\overline{\Pi}_{\mathcal{H}}(u) = \int_u^\infty \overline{\Pi}_S(z)\, dz, \qquad u > 0. \tag{6.2}$$

This implies that the integral in (6.2) and so $\mathbb{E}[X_1]$ is finite.

We will first investigate, which business line is most likely to cause ruin. Recall $P^k$ from (3.4) and the representation of their tail integrals in Corollary 5.5.

Invoking our quintuple law, we obtain precise asymptotic results on the ruin event. For this result, we need that $\overline{\Pi}_{\mathcal{H}}$ is subexponential; we write $\overline{\Pi}_{\mathcal{H}} \in \mathcal{S}$. If $\Pi_S$ is finite and has infinite support, this is implied by the d.f. of the increment $S_1 = S_1^1 + S_1^2$ to belong to the class $\mathcal{S}^*$ as introduced in Klüppelberg [16].

Next, we recall that subexponential distributions or d.f.s in $\mathcal{S}^*$ can belong to two different maximum domains of attraction $\mathrm{MDA}(\Phi_\alpha)$ for some $\alpha > 0$ or $\mathrm{MDA}(\Lambda)$. The first class covers the regular variation case, which has been investigated in [9]. From Theorem 3.1 of that paper, we know that $S_1 = S_1^1 + S_1^2$ is regularly varying, provided that the marginals are regularly varying in combination with a homogeneous Lévy copula. The second class contains subexponentials with lighter tails like lognormal or heavy-tailed Weibull distributions. For details see [11], Chapter 3.

THEOREM 6.1. *Suppose that $(X^1, X^2)$ is a spectrally positive Lévy process such that $X = X^1 + X^2$ satisfies the net profit condition (6.1). Assume that either:*

(i) $\overline{\Pi}_S \in \mathrm{RV}_{-\alpha}$ for $\alpha > 1$ or
(ii) $\overline{\Pi}_S \in \mathrm{MDA}(\Lambda) \cap \mathcal{S}$ and $\overline{\Pi}_{\mathcal{H}} \in \mathcal{S}$.

*Then the ruin probability is subexponential, i.e., $\mathbb{P}(\tau^+_\cdot < \infty) \in \mathcal{S}$.*

*[In case (i), we have $\mathbb{P}(\tau^+_\cdot < \infty) \in \mathrm{RV}_{-\alpha+1}$.]*

*Let $a(x) \sim \int_x^\infty \overline{\Pi}_S(z)\, dz / \overline{\Pi}_S(x)$ as $x \to \infty$ and suppose that the Lévy copula satisfies (5.9). Then for $k = 1, 2$ and $u, v > 0$, we have*

$$\lim_{x \to \infty} \mathbb{P}\left(\frac{X_{\tau_x^+} - x}{a(x)} > u, \frac{-X_{\tau_x^+-}}{a(x)} > v, \Delta X_{\tau_x^+} = \Delta P^k_{\tau_x^+} \,\Big|\, \tau_x^+ < \infty\right) = 0, \tag{6.3}$$

$$\lim_{x \to \infty} \mathbb{P}(\Delta X_{\tau_x^+} = \Delta P^k_{\tau_x^+} \mid \tau_x^+ < \infty) = 0 \tag{6.4}$$

*and*

$$\lim_{x \to \infty} \mathbb{P}\left(\frac{X_{\tau_x^+} - x}{a(x)} > u, \frac{-X_{\tau_x^+-}}{a(x)} > v, \Delta X_{\tau_x^+} = \Delta P^3_{\tau_x^+} \,\Big|\, \tau_x^+ < \infty\right)$$
$$= \mathrm{GPD}(u + v),$$

(6.5) $$\lim_{x \to \infty} \mathbb{P}(\Delta X_{\tau_x^+} = \Delta P^3_{\tau_x^+} \mid \tau_x^+ < \infty) = 1.$$

In case (i), $\text{GPD}(u+v) = (1+\frac{u+v}{\alpha})^{-\alpha}$ and $a(x) \sim x/\alpha$; in case (ii), $\text{GPD}(u+v) = e^{-(u+v)}$.

REMARK 6.2. (i) Theorem 6.1 generalizes the CPP situation in [5], Corollary 3.6, where the ruin probability was calculated for Pareto distributed jump sizes and a Clayton–Lévy copula.

(ii) By [18], Remark 4.3(iii), ruin can asymptotically occur for subexponential $\overline{\Pi}_{\mathcal{H}}$ only by a jump. In the situation of Theorem 6.1, relation (6.5) means that ruin occurs asymptotically only by a common jump, i.e., a claim that applies to both business lines at the same time.

The next result considers the barrier $x = 0$.

COROLLARY 6.3. Suppose that $(X^1, X^2)$ is a spectrally positive Lévy process such that 0 is irregular for $(0, \infty)$ with respect to $X = X^1 + X^2$.

(a) If the dependence is modeled by a Clayton–Lévy copula $\widehat{C}_\theta$, then
$$\lim_{\theta \to \infty} \mathbb{P}(\Delta X_{\tau_0^+} = \Delta P^k_{\tau_0^+} \mid \tau_0^+ < \infty) = \begin{cases} 0, & \text{for } k = 1, 2, \\ 1, & \text{for } k = 3. \end{cases}$$

(b) If the dependence is modeled by the nonhomogeneous Lévy copula $\widehat{C}_\zeta$, then
$$\lim_{\zeta \to \infty} \mathbb{P}(\Delta X_{\tau_0^+} = \Delta P^k_{\tau_0^+} \mid \tau_0^+ < \infty) = \begin{cases} \dfrac{\mu_{S^k}}{\mu_S}, & \text{for } k = 1, 2, \\ 0, & \text{for } k = 3. \end{cases}$$

## 7. Proofs.

PROOF OF THEOREM 3.1. *Case* 1. $S^1 + S^2$ is of bounded variation.

Let $m, k, f, g$ and $h$ be positive continuous functions with compact support satisfying $f(0) = g(0) = h(0) = 0$. The condition $f(0) = g(0) = h(0) = 0$ is to exclude the case of first passage by creeping, i.e., the event $\{X_{\tau_x^+} = x\}$, because we consider only the case, when the overshoot $X_{\tau_x^+} - x$ is a.s. positive. Since $S^1 + S^2$ is of bounded variation, we decompose it as in (3.3) into

$$S^1 + S^2 = P^1 + P^2 + P^3 + P^4 + P^5 + S^{1,-} + S^{2,-} + S^{--}.$$

Let $J_{S^1+S^2}$ denote the Poisson random measure associated with the jumps of $S^1 + S^2$. Then

$$J_{S^1+S^2} = J_{P^1} + J_{P^2} + J_{P^3} + J_{P^4} + J_{P^5} + J_{S^{1,-}} + J_{S^{2,-}} + J_{S^{--}},$$

where $J_{P^k}$ denotes the Poisson random measure associated with the jumps of $P^k$ given in (3.3). As $P^1, P^2, P^3, P^4, P^5, S^{1,-}, S^{2,-}$ and $S^{--}$ are independent, $J_{P^1}, J_{P^2}, J_{P^3}, J_{P^4}, J_{P^5}, J_{S^{1,-}}, J_{S^{2,-}}$ and $J_{S^{--}}$ have disjoint support with probability one. $J_{P^k}$ has intensity measure $\Pi_{P^k}(dx)\,dt$ and analogously to Step 1 of the proof of Theorem 3 in [8] we obtain for $k=1,\ldots,5$,

$$\int_{u>0,y\in[0,x],v\geq y,s\geq 0,t\geq 0} m(t)k(s)f(u)g(v)h(y)$$
$$\times \mathbb{P}(\tau_x^+ - \overline{G}_{\tau_x^+-} \in dt, \overline{G}_{\tau_x^+-} \in ds, X_{\tau_x^+} - x \in du,$$
$$x - X_{\tau_x^+-} \in dv, x - \overline{X}_{\tau_x^+-} \in dy, \Delta X_{\tau_x^+} = \Delta P^k_{\tau_x^+})$$
$$= \int_{y\in[0,x]}\int_{s\in[0,\infty)}\int_{v\in[y,\infty)}\int_{t\in[0,\infty)} m(t)k(s)h(y)g(v)$$
$$\times \int_{(0,\infty)} f(u)\Pi_{P^k}(du+v)\widehat{\mathcal{U}}(dt,dv-y)\mathcal{U}(ds,x-dy).$$

*Case* 2. $S^1 + S^2$ is of unbounded variation.

We start with the truncated process $(S^{1,\varepsilon}, S^{2,\varepsilon})$ for $\varepsilon > 0$ as given in (3.2). Then $S^{1,\varepsilon} + S^{2,\varepsilon}$ is of bounded variation and we can decompose its sample paths according to its jump behavior like in (3.3). The unbounded variation of $X$ implies that 0 is regular for $(0,\infty)$ and $(-\infty,0)$; see [19], Theorem 6.5(i). Therefore, $\mathcal{U}(\{0\},\{0\}) = \widehat{\mathcal{U}}(\{0\},\{0\}) = 0$. We apply the result of Case 1 above dropping the point 0 from integration, which yields

$$\int_{u>0,y\in[0,x],v\geq y,u+v>\varepsilon,s\geq 0,t\geq 0} m(t)k(s)f(u)g(v)h(y)$$
$$\times \mathbb{P}(\tau_x^+ - \overline{G}_{\tau_x^+-} \in dt, \overline{G}_{\tau_x^+-} \in ds,$$
$$X_{\tau_x^+} - x \in du, x - X_{\tau_x^+-} \in dv,$$
$$x - \overline{X}_{\tau_x^+-} \in dy, \Delta X_{\tau_x^+} = \Delta P^{k,\varepsilon}_{\tau_x^+})$$
$$= \int_{\phi\in(0,\infty)}\int_{t\in(0,\infty)}\int_{\xi\in(0,x]}\int_{s\in(0,\infty)} m(t)k(s)h(x-\xi)g(x+\phi-\xi)$$
$$\times \int_{(x+\phi-\xi,\infty)} f(\eta-(x+\phi-\xi))\Pi_{P^{k,\varepsilon}}(d\eta)\mathcal{U}(ds,d\xi)\widehat{\mathcal{U}}(dt,d\phi)$$
$$= \int_{y\in[0,x)}\int_{s\in(0,\infty)}\int_{v\in(y,\infty)}\int_{t\in(0,\infty)} m(t)k(s)h(y)g(v)$$
$$\times \int_{(v,\infty)} f(\eta-v)\Pi_{P^{k,\varepsilon}}(d\eta)\widehat{\mathcal{U}}(dt,dv-y)\mathcal{U}(ds,x-dy).$$

For $\varepsilon \downarrow 0$, the processes $P^{k,\varepsilon}$ converge a.s. to $P^k$ and, hence, in distribution. Moreover, for all $v > 0$ the function $\widetilde{f}(\eta) := f(\eta-v)1_{\{\eta>v\}}$ is bounded and continuous and vanishes on $[0, v]$. By [21], Theorem 8.7, for all $v > 0$

$$\lim_{\varepsilon \downarrow 0} \int_{(v,\infty)} f(\eta - v) \Pi_{P^{k,\varepsilon}}(d\eta) = \int_{(0,\infty)} f(u) \Pi_{P^k}(du + v)$$

and (3.9) follows. $\square$

PROOF OF COROLLARY 3.2. Since $X$ is not a CPP, its maxima are obtained at unique times, so $\overline{G}_{\tau_0^+ -} = \sup\{s < \tau_0^+ : X_s = 0\} = 0$ a.s. The proof for case (I) is analogous to Case 1 of the proof of Theorem 3.1. $\square$

PROOF OF THEOREM 4.3. To calculate $\mathcal{U}$ and $\widehat{\mathcal{U}}$ in Theorem 3.1 explicitly, we have to specify the local time at maximum and at minimum such that the normalization condition (3.8) is satisfied. Since $X$ is spectrally positive, we choose the local time at the minimum as

$$\widehat{L}_t = -\overline{X}_t = c \int_0^t 1_{\{\underline{X}_s = X_s\}} ds,$$

where $\underline{X}$ is defined in (1.2). The unkilled descending ladder process is for $t \geq 0$,

(7.1) $\qquad (\widehat{L}_t^{-1}, \widehat{H}_t) = (\inf\{s > 0 : \underline{X}_s < -t\}, \widehat{X}_{\widehat{L}_t^{-1}}) = (\tau_{-t}^-, t).$

Thus, by (3.7), we obtain

$$\widehat{\mathcal{U}}(ds, dx) = \int_0^\infty \mathbb{P}(\widehat{L}_t^{-1} \in ds, \widehat{H}_t \in dx) \, dt = \mathbb{P}(\widehat{L}_x^{-1} \in ds) \, dx$$
(7.2)
$$= \mathbb{P}(\tau_{-x}^- \in ds) \, dx,$$

$\widehat{\mathcal{U}}([0,\infty), dx) = dx$ and $\widehat{\kappa}(0, \beta) = \beta$, see (3.5). When the normalization condition (3.8) is satisfied, by Vigon [22], Proposition 3.3, it follows for $z > 0$,

(7.3) $\qquad \overline{\Pi}_{\mathcal{H}}(z) = \int_z^\infty \overline{\Pi}_S(x) \widehat{\mathcal{U}}([0,\infty), dx) = \int_z^\infty \overline{\Pi}_S(x) \, dx.$

Due to the irregularity of 0 for $[0, \infty)$ and following [19], Theorem 6.7(ii), we choose the local time at the maximum as

$$L_t = \sum_{k=0}^{n_t} e_\zeta^{(k)} \qquad \text{with } n_t = \#\{0 < s \leq t : \overline{X}_s = X_s\}$$

for an arbitrary parameter $\zeta > 0$ and i.i.d. $e_\zeta^{(k)} \stackrel{d}{=} \text{expo}(\zeta)$. Further, due to condition (4.5), the ascending ladder process is killed, i.e., there is a bivariate CPP $(\mathcal{L}^{-1}, \mathcal{H})$ with jump intensity $\zeta$ and $q > 0$ such that

$$\{(L_t^{-1}, H_t) : t < L_\infty\} \stackrel{d}{=} \{(\mathcal{L}_t^{-1}, \mathcal{H}_t) : t < e_q\}$$

and $(L_t^{-1}, H_t) = (\infty, \infty)$ for $t \geq L_\infty \stackrel{\mathrm{d}}{=} e_q$. $\mathcal{H}$ is a CPP with intensity $\zeta$ and with (7.3) the normalization condition (3.8) is satisfied if and only if $\zeta = \Pi_\mathcal{H}(\mathbb{R}) = \mu_S$. From (3.5), we obtain

$$\kappa(0, -i\theta) = q + \int_0^\infty (1 - e^{i\theta x}) \Pi_\mathcal{H}(dx)$$

and with (4.4) and the Wiener–Hopf factorization, see [19], equation (6.21), it becomes

$$\kappa(0, -i\theta) = k' \frac{\Psi_X(\theta)}{\widehat{\kappa}(0, i\theta)} = \frac{k'}{i\theta}\left(ic\theta + \int_0^\infty (1 - e^{i\theta x})\Pi_S(dx)\right).$$

Since $\mathcal{H}$ is of bounded variation and $\lim_{x \downarrow 0} x\overline{\Pi}_S(x) = 0$ by (7.3), partial integration results in

$$\kappa(0, -i\theta) = k'\left((c - \mu_S) + \int_0^\infty (1 - e^{i\theta x})\overline{\Pi}_S(x)\,dx\right).$$

From (7.3), we conclude $k' = 1$ and $q = c - \mu_S$. Since $e_q \stackrel{\mathrm{d}}{=} \operatorname{expo}(q)$ with $q = c - \mu_S$ is independent of $(\mathcal{L}^{-1}, \mathcal{H})$ and $\mathcal{N}_t = \#\{0 < s \leq t : \Delta\mathcal{H}_t \neq 0\}$ is a Poisson process with intensity $\zeta = \mu_S$, we get with $F_{\mathcal{L}^{-1}\mathcal{H}} = \frac{1}{\mu_S}\Pi_{\mathcal{L}^{-1}\mathcal{H}}$ for $s \geq 0, x \geq 0$

(7.4)
$$\mathcal{U}(ds, dx) = \int_0^\infty \mathbb{P}(t < e_q, \mathcal{L}_t^{-1} \in ds, \mathcal{H}_t \in dx)\,dt$$
$$= \frac{1}{c}\sum_{n=0}^\infty \left(\frac{\mu_S}{c}\right)^n F_{\mathcal{L}^{-1}\mathcal{H}}^{n*}(ds, dx).$$

Finally, from the quintuple law (3.1) with (7.2) and (7.4), we obtain for $u > 0$, $y \in [0, x]$, $v \geq y$, $s \geq 0$, $t \geq 0$ and for $k = 1, 2, 3$,

$$\mathbb{P}(\tau_x^+ - \overline{G}_{\tau_x^+-} \in dt, \overline{G}_{\tau_x^+-} \in ds, X_{\tau_x^+} - x \in du, x - X_{\tau_x^+-} \in dv,$$
$$x - \overline{X}_{\tau_x^+-} \in dy, \Delta X_{\tau_x^+} = \Delta P_{\tau_x^+}^k)$$
$$= \Pi_{P^k}(du + v)\widehat{\mathcal{U}}(dt, dv - y)\mathcal{U}(ds, x - dy)$$
$$= \Pi_{P^k}(du + v)\mathbb{P}(\tau_{-(v-y)}^- \in dt)\mathbf{1}_{\{v - y \geq 0\}}\,dv$$
$$\times \frac{1}{c}\sum_{n=0}^\infty \left(\frac{\mu_S}{c}\right)^n F_{\mathcal{L}^{-1}\mathcal{H}}^{n*}(ds, x - dy).$$

According to [8], Corollary 6, we have

$$\Pi_{\mathcal{L}^{-1}\mathcal{H}}(dt, dh) = \int_{[0,\infty)} \Pi_S(dh + \theta)\widehat{\mathcal{U}}(dt, d\theta)$$

and with (7.2) and the normalization condition expression (4.7) results. $\square$

PROOF OF COROLLARY 4.4. Integrating out time in (4.6) yields (4.8). Relation (4.11) follows from (7.3) with the normalization condition and the decomposition (3.3). The identity (4.9) results from (4.8) by integrating out $u$, $v$ and $y$. By integrating out $s$, $t$ and $v$ in the quintuple law in [8], Theorem 3, we obtain

$$\mathbb{P}(X_{\tau_x^+} - x \in du, x - \overline{X}_{\tau_x^+} \in dy) = \mathcal{U}(x - dy)\Pi_{\mathcal{H}}(du + y)$$

and with (7.4) we obtain

$$\mathbb{P}(\tau_x^+ < \infty) = \int_0^x \overline{\Pi}_{\mathcal{H}}(y)\mathcal{U}(x - dy)$$

$$= \sum_{n=0}^{\infty} \left(\frac{\mu_S}{c}\right)^{n+1} \int_0^x \overline{F}_{\mathcal{H}}(y) F_{\mathcal{H}}^{n*}(x - dy)$$

$$= \sum_{n=0}^{\infty} \left(\frac{\mu_S}{c}\right)^{n+1} (F_{\mathcal{H}}^{n*}(x) - F_{\mathcal{H}}^{(n+1)*}(x))$$

$$= \left(1 - \frac{\mu_S}{c}\right) \sum_{n=1}^{\infty} \left(\frac{\mu_S}{c}\right)^n \overline{F_{\mathcal{H}}^{n*}}(x). \quad \square$$

PROOF OF COROLLARY 4.6. For the barrier $x = 0$, we obtain under the normalization condition (3.8) with (3.10) and (7.4)

$$\mathbb{P}(\tau_0^+ \in dt, X_{\tau_0^+} \in du, -X_{\tau_0^+-} \in dv, \Delta X_{\tau_0^+} = \Delta P_{\tau_0^+}^k)$$

$$= \Pi_{P^k}(du + v)\widehat{\mathcal{U}}(dt, dv)\mathcal{U}(\{0\}, \{0\})$$

$$= \Pi_{P^k}(du + v)\mathbb{P}(\tau_{-v}^- \in dt)\,dv\frac{1}{c}.$$

The identities (4.13) and (4.14) follow from (4.12) by integrating out $t$, $u$ and $v$. Equation (4.15) results by summing up (4.14) for $k = 1, 2, 3$. $\square$

PROOF OF THEOREM 5.4. The tail integral of $(S^1, S^2)$ is given by

$$\overline{\Pi}(x_1, x_2) = \widehat{C}(\overline{\Pi}_{S^1}(x_1), \overline{\Pi}_{S^2}(x_2)), \qquad x_1, x_2 \in \mathbb{R} \setminus \{0\}.$$

So, we get for $z > 0$,

$$\overline{\Pi}_{P^1}(z) = \overline{\Pi}_{S^1}(z) - \lim_{y\downarrow 0} \widehat{C}(\overline{\Pi}_{S^1}(z), \overline{\Pi}_{S^2}(y)) + \lim_{y\uparrow 0} \widehat{C}(\overline{\Pi}_{S^1}(z), \overline{\Pi}_{S^2}(y)).$$

Analogous calculations give $\overline{\Pi}_{P^2}$.

For the common jump measures, we obtain for $z > 0$

(7.5) $\quad \overline{\Pi}_{P^3}(z) = \Pi(\{(x, y) \in (0, \infty) \times (0, \infty) : x + y > z\}),$

(7.6) $\quad \overline{\Pi}_{P^4}(z) = \Pi(\{(x, y) \in (0, \infty) \times (-\infty, 0) : x + y > z\}),$

(7.7) $\quad \overline{\Pi}_{P^5}(z) = \Pi(\{(x, y) \in (-\infty, 0) \times (0, \infty) : x + y > z\}).$

Since $\widehat{C}$ is twice continuously differentiable, we have (cf. [6], page 148)

$$(7.8) \qquad \Pi(dx, dy) = \frac{\partial^2 \widehat{C}(u,v)}{\partial u\, \partial v}\bigg|_{u=\overline{\Pi}_{S^1}(x), v=\overline{\Pi}_{S^2}(y)} \Pi_{S^1}(dx)\Pi_{S^2}(dy).$$

So the r.h.s. of (7.5) is given by

$$\int_0^\infty \int_{(z-x)\vee 0}^\infty \Pi(dx, dy)$$

$$= \int_0^\infty \int_{(z-x)\vee 0}^\infty \frac{\partial^2 \widehat{C}(u,v)}{\partial u\, \partial v}\bigg|_{u=\overline{\Pi}_{S^1}(x), v=\overline{\Pi}_{S^2}(y)} \Pi_{S^2}(dy)\Pi_{S^1}(dx)$$

$$= \int_0^\infty \frac{\partial \widehat{C}(u,v)}{\partial u}\bigg|_{u=\overline{\Pi}_{S^1}(x), v=\overline{\Pi}_{S^2}((z-x)\vee 0)} \Pi_{S^1}(dx),$$

since $\widehat{C}(u,0) = 0$ for all $u \in \mathbb{R}$. The r.h.s. of (7.6) and (7.7) are calculated analogously. $\square$

PROOF OF LEMMA 5.6. With Corollary 5.5, we get by homogeneity,

$$\frac{\overline{\Pi}_{P^1}(z)}{\overline{\Pi}_{S^1}(z)} = 1 - \lim_{y\downarrow 0} \widehat{C}\left(1, \frac{\overline{\Pi}_{S^2}(y)}{\overline{\Pi}_{S^1}(z)}\right), \qquad z \geq 0,$$

which is equal to 1 for all $z \geq 0$ in the case of independence. Otherwise, by left-continuity in $\infty$,

$$\lim_{z\to\infty} \frac{\overline{\Pi}_{P^1}(z)}{\overline{\Pi}_{S^1}(z)} = 1 - \lim_{z\to\infty}\lim_{y\downarrow 0} \widehat{C}\left(1, \frac{\overline{\Pi}_{S^2}(y)}{\overline{\Pi}_{S^1}(z)}\right) = 1 - \widehat{C}(1,\infty) = 0.$$

The proof for $P^2$ is analogous. $\square$

PROOF OF THEOREM 5.7. (1) Equation (5.10) holds by

$$\lambda = \Pi_S((0,\infty)) = \Pi([0,\infty)^2)$$

$$= \Pi_{S^1}((0,\infty)) + \Pi_{S^2}((0,\infty)) - \Pi((0,\infty)^2) = \lambda_1 + \lambda_2 - \lim_{x\downarrow 0}\overline{\Pi}(x,x)$$

and equation (5.11) results from (4.2) with Corollary 5.5.
(2) Equation (4.11) and Corollary 5.5 yield relation (5.12). $\square$

PROOF OF THEOREM 6.1. From [18], Lemma 3.5, we have for $\overline{\Pi}_{\mathcal{H}} \in \mathcal{S}$ the relation

$$(7.9) \qquad \lim_{x\to\infty} \frac{\mathbb{P}(\tau_x^+ < \infty)}{\overline{\Pi}_{\mathcal{H}}(x)} = \mathcal{U}([0,\infty)) = \frac{1}{|\mathbb{E}[X_1]|}.$$

In case (i) the assumption $\overline{\Pi}_S \in \mathrm{RV}_{-\alpha}$ and Karamata's theorem ([3], Theorem 1.5.11(ii)) to (6.2) yields $\overline{\Pi}_{\mathcal{H}} \in \mathrm{RV}_{-\alpha+1} \subset \mathcal{S}$. So, the first assertion results.

From Theorem 3.1, it follows for $u^*, v^* > 0$

$$\mathbb{P}(X_{\tau_x^+} - x > u^*, x - X_{\tau_x^+-} > v^*, \Delta X_{\tau_x^+} = \Delta P^1_{\tau_x^+})$$
$$= \int_{y \in [0,x]} \overline{\mu}_1(u^* + (v^* \vee y))\mathcal{U}(x - dy),$$

where $\overline{\mu}_1(z) := \int_z^\infty \overline{\Pi}_{P^1}(s)\, ds$. For $u, v > 0$, defining $u^* := a(x)u, v^* := x + a(x)v$, we have

$$\mathbb{P}\left(\frac{X_{\tau_x^+} - x}{a(x)} > u, \frac{-X_{\tau_x^+-}}{a(x)} > v, \Delta X_{\tau_x^+} = \Delta P^1_{\tau_x^+}\right)$$
$$= \mathcal{U}([0, x))\overline{\mu}_1(x + a(x)(u + v)).$$

With (7.9), we obtain

(7.10)
$$\lim_{x \to \infty} \mathbb{P}\left(\frac{X_{\tau_x^+} - x}{a(x)} > u, \frac{-X_{\tau_x^+-}}{a(x)} > v, \Delta X_{\tau_x^+} = \Delta P^1_{\tau_x^+} \,\Big|\, \tau_x^+ < \infty\right)$$
$$= \lim_{x \to \infty} \frac{\mathcal{U}([0,x))\overline{\mu}_1(x + a(x)(u+v))}{\mathcal{U}([0,\infty))\overline{\Pi}_{\mathcal{H}}(x)} \leq \lim_{x \to \infty} \frac{\mathcal{U}([0,x))\overline{\mu}_1(x)}{\mathcal{U}([0,\infty))\overline{\Pi}_{\mathcal{H}}(x)}.$$

Now, recall that

(7.11)
$$\frac{\overline{\mu}_1(x)}{\overline{\Pi}_{\mathcal{H}}(x)} = \frac{\int_x^\infty \overline{\Pi}_{P^1}(s)\, ds}{\int_x^\infty \overline{\Pi}_S(s)\, ds}.$$

Since by Lemma 5.6,

$$\frac{\overline{\Pi}_{P^1}(x)}{\overline{\Pi}_S(x)} \leq \frac{\overline{\Pi}_{P^1}(x)}{\overline{\Pi}_{S^1}(x)} \to 0 \qquad \text{as } x \to \infty$$

holds, the right-hand side of (7.11) tends to 0 as $x \to \infty$ by l'Hospital's lemma, and hence, the right-hand bound of (7.10) is 0. The analogous result holds for $P^2$. This implies (6.3) and (6.4). By [17], Theorems 1 and 2, we get for $u, v > 0$

$$\lim_{x \to \infty} \mathbb{P}\left(\frac{X_{\tau_x^+} - x}{a(x)} > u, \frac{-X_{\tau_x^+-}}{a(x)} > v \,\Big|\, \tau_x^+ < \infty\right) = \mathrm{GPD}(u + v).$$

Therefore, the last two relations result. $\square$

PROOF OF COROLLARY 6.3. If 0 is irregular for $(0, \infty)$, we get with Corollary 3.2

$$\mathbb{P}(\Delta X_{\tau_0^+} = \Delta P^k_{\tau_0^+}) = \int_0^\infty \overline{\Pi}_{P^k}(z)\, dz\, \mathcal{U}(\{0\})$$

and $\mathbb{P}(\tau_0^+ < \infty) = \mu_S \mathcal{U}(\{0\})$ where $S$ denotes the jump part of $X$ and $\mu_S = \mathbb{E}[S_1]$. Note that $\mathcal{U}(\{0\}) > 0$, if 0 is irregular for $(0, \infty)$. So

$$\mathbb{P}(\Delta X_{\tau_0^+} = \Delta P_{\tau_0^+}^k \mid \tau_0^+ < \infty) = \frac{\mu_{P^k}}{\mu_S},$$

where $\mu_{P^k} = \int_0^\infty \overline{\Pi}_{P^k}(z)\, dz$. From Section 5.3.3, we know that $\lim_{\theta \to \infty} \overline{\Pi}_{P^k}(z) = 0$ for $k = 1, 2$. Furthermore, increasing the dependence parameter $\zeta$ of the nonhomogeneous Lévy copula lowers the tail integral of the common jump, i.e., $\lim_{\zeta \to \infty} \overline{\Pi}_{P^3}(z) = 0$ and $\lim_{\zeta \to \infty} \overline{\Pi}_{P^k}(z) = \overline{\Pi}_{S^k}$ for $k = 1, 2$. Thus, the assertions follow. $\square$

**Acknowledgment.** We take pleasure in thanking Andreas Kyprianou for various interesting discussions.

## REFERENCES

[1] BERTOIN, J. (1996). *Lévy Processes. Cambridge Tracts in Mathematics* **121**. Cambridge Univ. Press, Cambridge. MR1406564
[2] BILLINGSLEY, P. (1979). *Probability and Measure.* Wiley, New York. MR534323
[3] BINGHAM, N. H., GOLDIE, C. M. and TEUGELS, J. L. (1987). *Regular Variation. Encyclopedia of Mathematics and Its Applications* **27**. Cambridge Univ. Press, Cambridge. MR898871
[4] BÖCKER, K. and KLÜPPELBERG, C. (2006). Multivariate models for operational risk. *Quant. Finance.* To appear.
[5] BREGMAN, Y. and KLÜPPELBERG, C. (2005). Ruin estimation in multivariate models with Clayton dependence structure. *Scand. Actuar. J.* **6** 462–480. MR2202587
[6] CONT, R. and TANKOV, P. (2004). *Financial Modelling with Jump Processes.* Chapman & Hall/CRC, Boca Raton, FL. MR2042661
[7] DONEY, R. A. (2007). *Fluctuation Theory for Lévy Processes. Lecture Notes in Math.* **1897**. Springer, Berlin. MR2320889
[8] DONEY, R. A. and KYPRIANOU, A. E. (2006). Overshoots and undershoots of Lévy processes. *Ann. Appl. Probab.* **16** 91–106. MR2209337
[9] EDER, I. and KLÜPPELBERG, C. (2007). Pareto Lévy measures and multivariate regular variation. To appear.
[10] ESMAEILI, H. and KLÜPPELBERG, C. (2008). Parameter estimation of a bivariate compound Poisson process. To appear.
[11] EMBRECHTS, P., KLÜPPELBERG, C. and MIKOSCH, T. (1997). *Modelling Extremal Events: For Insurance and Finance. Applications of Mathematics (New York)* **33**. Springer, Berlin. MR1458613
[12] HUZAK, M., PERMAN, M., ŠIKIĆ, H. and VONDRAČEK, Z. (2004). Ruin probabilities and decompositions for general perturbed risk processes. *Ann. Appl. Probab.* **14** 1378–1397. MR2071427
[13] HUZAK, M., PERMAN, M., ŠIKIĆ, H. and VONDRAČEK, Z. (2004). Ruin probabilities for competing claim processes. *J. Appl. Probab.* **41** 679–690. MR2074816
[14] JOE, H. (1997). *Multivariate Models and Dependence Concepts. Monographs on Statistics and Applied Probability* **73**. Chapman & Hall, London. MR1462613
[15] KALLSEN, J. and TANKOV, P. (2006). Characterization of dependence of multidimensional Lévy processes using Lévy copulas. *J. Multivariate Anal.* **97** 1551–1572. MR2275419

[16] KLÜPPELBERG, C. (1988). Subexponential distributions and integrated tails. *J. Appl. Probab.* **25** 132–141. MR929511
[17] KLÜPPELBERG, C. and KYPRIANOU, A. (2006). On extreme ruinous behaviour of Lévy insurance risk processes. *J. Appl. Probab.* **43** 1–5.
[18] KLÜPPELBERG, C., KYPRIANOU, A. E. and MALLER, R. A. (2004). Ruin probabilities and overshoots for general Lévy insurance risk processes. *Ann. Appl. Probab.* **14** 1766–1801. MR2099651
[19] KYPRIANOU, A. E. (2006). *Introductory Lectures on Fluctuations of Lévy Processes with Applications.* Springer, Berlin. MR2250061
[20] NELSEN, R. B. (2006). *An Introduction to Copulas*, 2nd ed. Springer, New York. MR2197664
[21] SATO, K. (1999). *Lévy Processes and Infinitely Divisible Distributions.* Cambridge Univ. Press, Cambridge.
[22] VIGON, V. (2002). Votre Lévy rampe-t-il? *J. London Math. Soc. (2)* **65** 243–256. MR1875147

CENTRE FOR MATHEMATICAL SCIENCES
AND INSTITUTE FOR ADVANCED STUDY
TECHNISCHE UNIVERSITÄT MÜNCHEN
BOLTZMANNSTRASSE 3
D-85747 GARCHING
GERMANY
E-MAIL: eder@ma.tum.de
cklu@ma.tum.de